\documentclass[graybox]{svmult}

\usepackage{type1cm}        
\usepackage{makeidx}         
\usepackage{graphicx}        
\usepackage{multicol}        
\usepackage[bottom]{footmisc}

\usepackage{newtxtext}       %
\usepackage{newtxmath}       

\makeindex             

\newcommand{\oN}{\mathbb{N}}
\newcommand{\oR}{\mathbb{R}}
\newcommand{\ignore}[1]{}

\begin{document}

\title*{A survey of semidefinite programming approaches to the generalized problem of moments and their error analysis}
\titlerunning{A survey of SDP approaches to the GPM}

\author{Etienne de Klerk        $\;\; \cdot \;\;$
        Monique Laurent}
\authorrunning{A survey of SDP approaches to the GPM} 
\institute{Etienne de Klerk \at Tilburg University and Delft University of Technology, \email{E.deKlerk@UvT.nl}
\and Monique Laurent \at CWI Amsterdam, and Tilburg University, \email{M.Laurent@CWI.nl}}

%
%
\maketitle

\abstract{The generalized problem of moments is a conic linear optimization problem over the convex cone of positive Borel measures
with given support. It has 
a large variety of applications, including global optimization of polynomials and rational functions, options pricing in finance, constructing quadrature schemes for numerical integration, and distributionally robust optimization.
A usual solution approach, due to J.B. Lasserre,  is to approximate the convex cone of positive Borel measures by  finite dimensional outer and inner conic approximations.
We will review some results on these approximations, with a special focus on the convergence rate of the hierarchies of upper and lower bounds for the general problem of moments that are obtained from these  inner and outer approximations.}

\section{Introduction}
\label{intro}

The classical problem of moments asks when a measure is determined by a set of specified moments  and  variants of this problem
were studied (in the univariate case) by leading 19th and early 20th century mathematicians, like Hamburger, Stieltjes, Chebyshev, Hausdorff, and Markov. We refer to \cite{Akh} for an early reference and to the
recent monograph  \cite{Sch} for a comprehensive treatment of the moment problem.

The generalized problem of moments asks to optimize a linear function over the set of finite, positive Borel measures that satisfy certain moment-type conditions.
More precisely, we consider continuous functions $f_0$ and $f_i$ ($i \in [m]$) where $[m] = \{1,\ldots,m\}$,
that are defined on a compact set $K \subset \mathbb{R}^n$.  The generalized problem
of moments (GPM) may now be defined as follows.\footnote{We only deal with the GPM in a restricted setting; more general versions
of the problem are studied in, e.g.,\, \cite{Shapiro2001}.}

\begin{trailer}{Generalized problem of moments (GPM)}
\begin{equation}
\label{prob:GPM}
val := \inf_{\mu \in \mathcal{M}(K)_+} \left\{ \int_K f_0(x)d\mu(x) \; : \; \int_K f_i(x)d\mu(x) = b_i \quad \forall i \in [m]\right\},
\end{equation}
where
\begin{itemize}
\item
$\mathcal{M}(K)_+$ denotes the convex cone of positive, finite, Borel measures (i.e.,\ Radon measures) supported on the set $K$;\footnote{Formally, we consider the
usual Borel $\sigma$-algebra, say $\mathcal{B}$,  on $\mathbb{R}^n$, i.e.,\ the smallest (or coarsest) $\sigma$-algebra that contains the open sets in $\mathbb{R}^n$.
A positive, finite Borel measure $\mu$ is a nonnegative-valued set function on $\mathcal{B}$, that is countably additive for disjoint sets in $\mathcal{B}$.
The support of $\mu$ is the set, denoted  $\text{Supp}(\mu)$, and  defined as the smallest closed set $S$ such that $\mu(\mathbb{R}^n\setminus S) = 0.$}
\item
The scalars $b_i\in\oR$ ($i \in [m]$) are given.
\end{itemize}
\end{trailer}

In this survey we will mostly consider the case where all $f_i$'s are polynomials, and will always assume $K \subseteq \mathbb{R}^n$ to be compact.
Moreover, for some of the results, we will  also assume that $K$ is a basic semi-algebraic set and we will sometimes further restrict to simple sets like a hypercube, simplex or sphere.

The generalized problem of moments has a rich history; see, e.g., \cite{Akh, Landau,Sch} and references therein and  \cite{Las18} for a recent overview of many of its applications. In the recent years modern optimization approaches have been investigated
in depth, in particular, by Lasserre (see \cite{Las2008}, the monograph \cite{lass-book} and further references therein).
Among others, there is a well-understood duality theory, and hierarchies of inner and outer approximations for the cone $\mathcal{M}(K)_+$ have been introduced that lead to converging upper and lower bounds for the problem (\ref{prob:GPM}).
In this survey we will present these hierarchies and  show how the corresponding bounds can be computed using semidefinite programming. Since several overviews are already available on general properties of these hierarchies (e.g., in \cite{lass-book,lass-book2,Lau09,Lau14}), our main focus here will be on recent results that describe their rate of convergence.
We will review in particular in more detail recent results  on the upper bounds arising from the inner approximations, and highlight some recent links made with orthogonal polynomials and cubature rules for integration.

\subsection{The dual problem of the GPM}

The GPM is an infinite-dimensional conic linear program, and therefore it has an associated dual problem.
Formally we introduce a duality (or pairing) between the following two 
vector spaces:
\begin{enumerate}
\item
the space $\mathcal{M}(K)$  of all signed, finite, Borel measures supported on $K$, 
\item
the space $\mathcal{C}(K)$ of
continuous functions on $K$, endowed with the supremum norm $\|\cdot\|_\infty$.
\end{enumerate}
The duality (pairing) in question is provided by the nondegenerate bilinear form $\langle \cdot,\cdot \rangle :\mathcal{C}(K) \times \mathcal{M}(K) \rightarrow \mathbb{R}$,
defined by
\[
\langle f, \mu \rangle = \int_K f(x)d\mu(x) \quad (f \in \mathcal{C}(K), \; \mu \in \mathcal{M}(K)).
\]
Thus the dual cone of $\mathcal{M}(K)_+$ w.r.t.\ this duality is the cone of continuous  functions that are nonnegative on $K$, and will
be denoted by $\mathcal{C}(K)_+ = \left(\mathcal{M}(K)_+ \right)^*$.

In our setting of compact $K \subset \mathbb{R}^n$, $\mathcal{M}(K)$ is also the  dual space of $\mathcal{C}(K)$, i.e.,\ $\mathcal{M}(K)$ may be associated with the
space of linear functionals defined on $\mathcal{C}(K)$. In particular, due to the  Riesz-Markov-Kakutani representation theorem (e.g. \cite[\S1.10]{Tao_epsilon}),
every linear functional on $\mathcal{C}(K)$ may be expressed as
\[
f \mapsto \langle f, \mu \rangle \quad \mbox{ for a suitable $\mu \in \mathcal{M}(K)$}.
\]
As a result, we have the weak$*$ topology on $\mathcal{M}(K)$ where the open sets are finite intersections of elementary sets
of the form
\[
\{ \mu \in \mathcal{M}(K) \; | \; \alpha <  \langle f, \mu \rangle < \beta\},
\]
for given  $\alpha,\beta \in \mathbb{R}$, and $f \in C(K)$, and the unions of such finite intersections.

A sequence $\{\mu_k\} \subset \mathcal{M}(K)$ converges in the weak$*$ topology, say $\mu_k \rightharpoonup \mu$,
if, and only if,
\begin{equation}
\label{eq:weak star convergence}
\lim_{k \rightarrow \infty} \langle f, \mu_k \rangle = \langle f, \mu \rangle \; \forall f \in \mathcal{C}(K).
\end{equation}
As a consequence of \eqref{eq:weak star convergence}, {the cone $\mathcal{M}(K)_+$ is closed} and the set of probability measures in $\mathcal{M}(K)$ is closed.

By Alaoglu's theorem, e.g.\ \cite[Theorem III(2.9)]{Barvinok_course in_convexity}, the following set {(i.e., the unit ball in $\mathcal{M}(K)$)} is compact
in the weak$*$ topology of $\mathcal{M}(K)$:
\begin{equation}
\label{unit ball in MK}
\left\{\mu \in \mathcal{M}(K) \; | \; |\langle f, \mu \rangle | \le 1 \; \forall f \in \mathcal{C}(K) \mbox{ with } \|f\|_\infty \le 1  \right\}.
\end{equation}
{Hence the set of probability measures in $\mathcal{M}(K)$ is compact, {since it is a closed subset of the compact set in (\ref{unit ball in MK})}, and thus 
 it provides a compact base in the weak$*$ topology for the cone $\mathcal{M}(K)_+$. This implies again that $\mathcal{M}(K)_+$ is closed in this topology (using Lemma 7.3 in \cite[Part IV]{Barvinok_course in_convexity}) and 
we will also use this fact  to analyze duality in the next section.}


 \begin{trailer}{Dual linear optimization problem of \eqref{prob:GPM}}
Using this duality setting, the dual conic linear program of (\ref{prob:GPM}) reads
\begin{eqnarray}
\notag
val^* &:=&\sup_{y \in \mathbb{R}^m} \left\{ \sum_{i \in [m]} b_i y_i \; : \; f_0 - \sum_{i \in I} y_if_i \in \mathcal{C}(K)_+ \right\}, \\
\label{dual:GPM}
       &=& \sup_{y \in \mathbb{R}^m} \left\{ \sum_{i \in [m]} b_i y_i \; : \; f_0(x) - \sum_{i \in I} y_if_i(x)  \ge 0 \; \forall x \in K \right\}.
\end{eqnarray}
\end{trailer}

By the duality theory of conic linear optimization, one has the following duality relations; see, e.g.,\ \cite[Section IV.7.2]{Barvinok_course in_convexity} or \cite[Appendix C]{lass-book}. 


\begin{theorem}\label{thm:sduality1}
Consider the GPM (\ref{prob:GPM}) and its dual (\ref{dual:GPM}). Assume (\ref{prob:GPM}) has a feasible solution.
One has $val \ge val^*$ (weak duality), with equality $val=val^*$ (strong duality) if the cone $\left\{ \left( \langle f_0,\mu\rangle , \langle f_1,\mu\rangle, \ldots \langle f_m,\mu\rangle\right) \; : \; \mu \in \mathcal{M}(K)_+  \right\}$
is a closed subset of $\mathbb{R}^{m+1}$.
 If, in addition, $val > -\infty$ then \eqref{prob:GPM}
 has an optimal solution.
\end{theorem}

We  mention another sufficient condition for strong duality, that is a consequence of Theorem \ref{thm:sduality1} in our setting. 

\begin{corollary}
 \label{thm:sduality2}
Assume   (\ref{prob:GPM}) has a feasible solution,  and there exist  
$z_0,z_1,\ldots,z_m\in\oR$  for which the function $\sum_{i=0}^mz_if_i$ is strictly positive on $K$ (i.e.,
$\sum_{i=0}^mz_i f_i(x)>0$ for all $x\in K$).
Then, $val=val^*$ holds and \eqref{prob:GPM} has an optimal solution.
\end{corollary}

Hence, if  in problem (\ref{prob:GPM}) we optimize over the probability measures (i.e., with $f_1\equiv 1$, $b_1=1$) 
 then the assumptions in Corollary \ref{thm:sduality2} are satisfied.

We indicate how Corollary \ref{thm:sduality2} can be derived from Theorem \ref{thm:sduality1}.
 Consider the linear map $L:\mathcal{M}(K)\to\oR^{m+1}$ defined by $L(\mu)=
(\langle f_0,\mu\rangle, \ldots,\langle f_m,\mu\rangle)$, which is continuous {w.r.t. the weak* topology on $\mathcal{M}(K)$.}
First we claim  $\text{Ker } L\cap \mathcal{M}(K)_+=\{0\}$. Indeed, assume $L(\mu)=0$ for some $\mu\in\mathcal{M}(K)_+$.
Setting $f=\sum_{i=0}^mz_if_i$, $L(\mu)=0$ implies  $\langle f,\mu \rangle =0$ and thus $\mu=0$ since $f$ is strictly positive on $K$.
Since the cone $\mathcal{M}(K)_+$ {has a compact convex base }
in the weak$*$ topology and the linear map $L$ is continuous, we can conclude that the image $L(\mathcal{M}(K)_+)$ is  closed  (using Lemma 7.3 in \cite[Part IV]{Barvinok_course in_convexity}). Now we can conclude using Theorem \ref{thm:sduality1}.

\subsection{Atomic solution of the GPM }

If the GPM has an optimal solution, then it has a finite atomic optimal solution, supported
on at most $m$ points (i.e.,\ the weighted sum of at most $m$ Dirac delta measures).
This  is a classical result in the theory of moments; see, e.g., \cite{rogosinski} (univariate case),  \cite{Kemperman} (which shows  an atomic measure with $m+1$ atoms using induction on $m$) and a modern exposition in \cite{Shapiro2001} (which shows an atomic measure with $m$ atoms).
The result may also
 be obtained as a consequence of the following, dimension-free version of the Carath\'eodory theorem.

\begin{theorem}[see, e.g.,\ Theorem 9.2 in Chapter III of \cite{Barvinok_course in_convexity}]\label{thm:Bar}
Let $S$ be a convex subset of a vector space such that, for every line $L$, the intersection
$S \cap L$ is a closed bounded interval. Then every extreme point of the intersection of $S$ with $m$ hyperplanes
can be expressed as a convex combination of at most $m+1$ extreme points of $S$.
\end{theorem}



 \begin{trailer}{Atomic solution of the (GPM)}
\begin{theorem} 
\label{thm:Rogosinsky}
If the GPM \eqref{prob:GPM} has an optimal solution then  it has one which is finite atomic with at most $m$ atoms, i.e., of the form
$
\mu^* = \sum_{\ell=1}^{m} w_\ell\delta_{x^{(\ell)}}
$
where $w_\ell \ge 0$, $x^{(\ell)} \in K$, and $\delta_{x^{(\ell)}}$ denotes the Dirac measure supported at $x^{(\ell)}$ ($\ell \in [m]$).
\end{theorem}
\end{trailer}



This result can be derived from Theorem \ref{thm:Bar} in the following way.
By assumption, the GPM has  an optimal solution $\mu^*$. Moreover, since it has one at an extreme point we may assume that  $\mu^*$ is an extreme point of the feasibility region $\mathcal{M}(K)_+\cap\cap_{i=1}^m H_i$ of the program (\ref{prob:GPM}), where $H_i$ is the hyperplane $\langle f_i,\mu\rangle =b_i$.
Then the following set $S = \{ \mu \in \mathcal{M}(K)_+ \; : \; \mu(K) =\mu^*(K)\}$ meets the condition of Theorem~\ref{thm:Bar}, since the 
set of probability measures in $\mathcal{M}(K)_+$ is compact in the weak$*$ topology, and any line in a topological vector space is closed (e.g.\ \cite[p. 111]{Barvinok_course in_convexity}).
Moreover, the extreme points of $S$ are precisely the scaled Dirac measures supported by  points in $K$ (see, e.g., Section III.8 in \cite{Barvinok_course in_convexity}).
In  addition,  $\mu^*$ is an extreme point of the set $S\cap\cap_{i=1}^m H_i$ and thus,
by Theorem~\ref{thm:Bar},  $\mu^*$ is a conic combination of $m+1$ Dirac measures supported at points
$x^{(\ell)}\in K$ for $\ell\in [m+1]$. Finally, as in \cite{Shapiro2001},  consider the LP
$$ \min \sum_{\ell=1}^{m+1} w_\ell f_0(x^{(\ell)}) \text{ s.t. } \ w_\ell\ge 0\ (\ell\in [m+1]), \sum_{\ell=1}^{m+1}w_\ell f_i(x^{(\ell)}) =b_i \ (i\in [m])$$
whose optimum value is $val$. Then an optimal solution attained at an extreme point provides an optimal solution of the GPM (\ref{prob:GPM}) which is atomic with at most $m$ atoms.




\subsection{GPM in terms of moments}

From now on we will assume  the functions $f_0,f_1,\ldots,f_m$ in the definition of the GPM (\ref{prob:GPM}) are all polynomials and the set $K$ is compact. Then the GPM may be reformulated
in terms of the moments of the variable measure $\mu$. To be precise, given a multi-index $\alpha = (\alpha_1,\ldots,\alpha_n) \in \mathbb{N}^n$
the moment of order $\alpha$ of a measure  $\mu \in \mathcal{M}(K)_+ $  is defined as
\[
m^\mu_\alpha(K) := \int_K x^\alpha d\mu(x).
\]
Here we set $x^\alpha=x_1^{\alpha_1}\cdots x_n^{\alpha_n}$.
 We may write the polynomials $f_0,f_1,\ldots,f_m$ in terms of the standard
monomial basis as:
\[
f_i(x) = \sum_{\alpha \in \mathbb{N}_d^n} f_{i,\alpha}x^\alpha \quad \forall i = 0,\ldots,m,
\]
where the $f_{i,\alpha} \in \mathbb{R}$ are the coefficients in the monomial basis, and we assume the maximum total degree of the polynomials
$f_0,f_1,\ldots,f_m$ to be at most $d$.

Throughout  we let $\oN^d_d=\{\alpha\in \oN^n: |\alpha|\le d\}$ denote the set of multi-indices, with
$|\alpha|=\sum_{i=1}^n\alpha_i$, and  $\oR[x]_d$ denotes the set of multivariate polynomials with degree at most $d$.

 \begin{trailer}{GPM in terms of moments}
 We may now rewrite the GPM \eqref{prob:GPM} in terms of moments:
\[
\inf_{\mu \in \mathcal{M}(K)_+} \left\{ \sum_{\alpha \in \mathbb{N}^n_d} f_{0,\alpha}m^\mu_\alpha(K) \; : \; \sum_{\alpha \in \mathbb{N}^n_d} f_{i,\alpha}m^\mu_\alpha(K) = b_i \; \forall i \in [m]\right\}.
\]
Here $d$ is the maximum  degree of the polynomials $f_0,f_1,\ldots,f_m$.
\end{trailer}

 Thus we  may consider the set of all possible truncated
moments sequences:
\[
 \left\{\left(m^\mu_\alpha(K)\right)_{\alpha \in \mathbb{N}_d^n} \; : \; \mu \in \mathcal{M}(K)_+ \right\},
\]
and  describe the inner and outer approximations for $\mathcal{M}(K)_+$ in terms of this set.

\subsection{Inner and outer approximations}\label{sec:approx}

We will consider two types of approximations of  the cone $\mathcal{M}(K)_+$, namely inner and outer conic approximations.

\paragraph{Inner approximations}
The underlying idea, due to Lasserre \cite{Las11}, is to consider a subset of measures $\mu$ in $\mathcal{M}(K)_+$ of the form
\[
d\mu = h \cdot d\mu_0,
\]
where $h$ is a polynomial sum-of-squares density function, and $\mu_0 \in \mathcal{M}(K)_+$ is a fixed reference measure with $\text{Supp}(\mu_0)=K$.

To obtain a finite dimensional subset of measures, we will limit the total degree of $h$ to some value $2r$ where $r \in \mathbb{N}$ is fixed.
The cone of sum-of-squares polynomials of total degree at most $2r$ will be denoted by $\Sigma_r$, hence
$$\Sigma_r=\left\{\sum_{i=1}^k p_i^2: k\in \oN, p_i\in \oR[x]_r, i\in [k]\right\}.$$
In this way one obtains the cones
\begin{equation}\label{inner approximation}
\mathcal{M}_{\mu_0}^r :=\left\{ \mu \in \mathcal{M}(K)_+ \; : \; d\mu = h \cdot d\mu_0, \; h \in \Sigma_r \right\}      \quad (r = 1,2, \ldots)
\end{equation}
which provide a hierarchy of inner approximations for the set $\mathcal{M}(K)_+$:
$$
\mathcal{M}^r_{\mu_0}\subseteq \mathcal{M}^{r+1}_{\mu_0} \subseteq \mathcal {M}(K)_+.
$$

\paragraph{Outer approximations}

The dual GPM (\ref{dual:GPM}) involves the nonnegativity constraint
\[
f_0(x) - \sum_{i=1}^m y_if_i(x)  \ge 0 \; \forall x \in K,
\]
which one may relax  to a sufficient condition that guarantees the nonnegativity
of the polynomial $f_0 - \sum_{i=1}^m y_if_i$ on $K$.
Lasserre \cite{lass-siopt-01} suggested to use the following sufficient condition
in the case when $K$ is a basic closed semi-algebraic set, i.e., when  we have a description of $K$ as the intersection of the level sets of polynomials $g_j$ ($j \in [k]$):
\[
K = \left\{x \in \mathbb{R}^n \; : \; g_j(x) \ge 0 \quad \forall j \in [k]\right\}.
\]
Namely, consider the condition
\[
f_0 - \sum_{i=1}^m y_if_i = \sigma_0 + \sum_{j=1}^k \sigma_j g_j,
\]
where each $\sigma_j$ is a sum-of-squares  polynomial and the degree of each term $\sigma_jg_j$ ($0\le j\le k$) is at most $2r$, so that the  degree of the right-hand-side polynomial is at
most $2r$. Here we set $g_0\equiv 1$ for notational convenience.
Thus we replace the cone $\mathcal{C}(K)_+$ by a cone of the type:
\begin{equation}\label{eq:quadratic module}
 \mathcal{Q}^r(g_1,\ldots,g_k) := \left\{ f  \; : \; f = \sigma_0 + \sum_{j=1}^k \sigma_j g_j, \; \sigma_j \in
  \Sigma_{r_j},\  j=0,1,\ldots,k\right\},
\end{equation}
where we set $r_j := r - \left\lceil\mbox{deg$(g_j)$}/2\right\rceil$ for all $ j \in \{0,\ldots,m\}$.

The cone $\mathcal{Q}^r(g_1,\ldots,g_k)$ is known as the \emph{truncated quadratic module} generated by the polynomials $g_1,\ldots,g_k$.
By definition, its dual cone consists of the signed measures $\mu$ supported on $K$ such that $\int_K f d\mu \ge 0$ for all $f \in \mathcal{Q}^r(g_1,\ldots,g_k)$:
\begin{equation}\label{dual:Qr}
(\mathcal{Q}^r(g_1,\ldots,g_k))^*=\left\{ \mu\in \mathcal{M}(K): \int_K f(x)d\mu(x)\ge 0 \quad \forall f\in \mathcal{Q}^r(g_1,\ldots,g_k)\right\}.
\end{equation}
This provides a hierarchy of outer approximations for the cone $\mathcal{M}(K)_+$:
$$
\mathcal{M}(K)_+\subseteq (\mathcal{Q}^{r+1}(g_1,\ldots,g_k))^* \subseteq (\mathcal{Q}^r(g_1,\ldots,g_k))^*.
$$
We will also briefly consider  the tighter outer approximations  for the cone $\mathcal{M}(K)_+$ obtained by replacing the
truncated quadratic module $\mathcal{Q}^r(g_1,\ldots,g_k)$ by the larger cone
$\mathcal{Q}^r\left(\prod_{j\in J}g_j: J\subseteq [k]\right)$,
thus the truncated quadratic module generated by all pairwise products of the $g_j$'s (also known as the pre-ordering generated by the $g_j$'s). Then we have
$$
\mathcal{M}(K)_+\subseteq \left(\mathcal{Q}^{r}\left(\prod_{j\in J}g_j: J\subseteq [k]\right)   \right )^* \subseteq (\mathcal{Q}^r(g_1,\ldots,g_k))^*.
$$

\ignore{
We will also briefly consider  the stronger bounds on $val$, denoted by $\overline{val}^{(r)}_{outer}$, that are obtained by replacing in the above definition of $val^{(r)}_{outer}$ the truncated quadratic module $\mathcal{Q}^r(g_1,\ldots,g_k)$ by the larger cone
$\mathcal{Q}^r\left(\prod_{j\in J}g_j: J\subseteq [k]\right)$,
thus the truncated quadratic module generated by all pairwise products of the $g_j$'s (also known as the pre-ordering generated by the $g_j$'s). Clearly we have
$$val^{(r)}_{outer} \le \overline{val}^{(r)}_{outer}\le val.$$
This new parameter can also be reformulated as a semidefinite program, analogous to the program (\ref{eq:upr1})-(\ref{eq:upr2}), but it now involves $2^k+1$ semidefinite constraints instead of $k+1$ such constraints in (\ref{eq:upr2}) and thus its practical implementation is feasible only for small values of $k$.
On the other hand, as we will see later in this section the bounds $\overline{val}^{(r)}_{outer}$ admit a much sharper error analysis than the bounds $val^{(r)}_{outer}$ for the case of polynomial optimization.
}

\section{Examples of GPM}
The GPM \eqref{prob:GPM} has many applications. Below we will list some examples that are directly relevant to this survey; additional
examples in control theory, options pricing in finance, and others, can be found in \cite{Las2008,lass-book,Las18}.

\begin{example}{Global minimization of polynomials on compact sets}
Consider the global optimization problem:
\begin{equation}
\label{eq:global opt pol}
val =  \min_{x \in K} p(x)
\end{equation}
where $p$ is  a polynomial and $K$ a compact set. 
This corresponds to
the GPM \eqref{prob:GPM} with $m=1$, $f_0 = p$, $f_1 = 1$ and $b_1 = 1$, i.e.:
\[
val = \min_{\mu \in \mathcal{M}(K)_+} \left\{ \int_K p(x)d\mu(x) \; : \; \int_K  d\mu(x) = 1\right\}.
\]
In the following sections we will focus on deriving error bounds for this problem when using the inner and outer approximations of
$\mathcal{M}(K)_+$.

\end{example}

\begin{example}{Global minimization of rational functions on compact sets}
We may generalize the previous example to rational objective functions.
In particular, we now consider the global optimization problem:
\begin{equation}\label{eq:rational global opt}
val =  \min_{x \in K} \frac{p(x)}{q(x)},
\end{equation}
where $p,q$ are polynomials such that $q(x) > 0$ $\forall$ $x \in K$, and $K \subseteq \mathbb{R}^n$ is  compact.

This problem has applications in many areas, including signal recovery \cite{Castella} and finding minimal energy configurations of point charges in a field with polynomial potential \cite{Thompson's problem}.

It is simple to see that we may reformulate this problem as the GPM with $m=1$ and $f_0 = p$,  $f_1 = q$, and $b_1 =1$, i.e.:
\[
val = \min_{\mu \in \mathcal{M}(K)_+} \left\{ \int_K p(x)d\mu(x) \; : \; \int_K q(x)d\mu(x) = 1\right\}.
\]
Indeed, one may readily verify that if $x^*$ is a global minimizer of the rational function $p(x)/q(x)$ over $K$ then
an optimal solution of the GPM is given by $\mu^* = \frac{1}{q(x^*)}\delta_{x^*}$.

\end{example}

\begin{example}{Polynomial cubature rules}
Positive cubature (also known as multivariate quadrature) rules for numerical integration of a function $f$ with respect to a measure $\mu_0$ over a set $K$ take the form
\[
\int_K f(x)d\mu_0(x) \approx \sum_{\ell=1}^N w_\ell f(x^{(\ell)}),
\]
where  the points $x^{(\ell)}\in K$  and the weights $w_\ell \ge 0$ ($\ell \in [N]$) are fixed. The points (also known as the nodes of the cubature rule) and weights are typically chosen so that the
approximation is exact for polynomials up to a certain degree, say $d$.

The problem of finding the points  $x^{(\ell)} \in K$ and weights $w_\ell$ ($\ell\in [N]$) giving a cubature rule exact at degree $d$ may then be written as the following GPM:
\[
val :=\inf_{\mu \in \mathcal{M}(K)_+} \left\{ \int_K 1 d\mu(x) \; : \; \int_K x^\alpha d\mu(x) = \int_K x^\alpha d\mu_0(x) \; \forall \alpha \in \mathbb{N}^n_d \right\}.
\]

The key observation is that, by Theorem \ref{thm:Rogosinsky}, this problem has an atomic solution supported on at most $N =  |\mathbb{N}^n_d| =  {n+d \choose d}$
points in $K$, say $\mu^*  = \sum_{\ell=1}^{N} w_\ell\delta_{x^{(\ell)}}$, and this yields the cubature  weights and points.
This result is known as Tchakaloff's theorem \cite{Tchakaloff}; see also \cite{BT,Trefethen cubature review}. (In fact,
our running assumption that $K$ is compact may be relaxed somewhat in Tchakaloff's theorem  --- see, e.g.\ \cite{Putinar cubature}.)

Here we have chosen the constant polynomial 1 as objective function so that the optimal value is $val= \mu_0(K)$. Other choices of objective functions are possible as discussed, e.g., in \cite{Ryu_Boyd_quadrature}.
The  GPM formulation of the cubature problem was used for the numerical calculation of cubature schemes for various sets $K$ in \cite{Ryu_Boyd_quadrature}.
\end{example}

\section{Semidefinite programming reformulations of the approximations}

The inner and outer approximations of the cone $\mathcal{M}(K)_+$ discussed in Section \ref{sec:approx}  lead to upper and lower bounds for the GPM (\ref{prob:GPM}), which may be reformulated as finite-dimensional, convex optimization problems,
namely semidefinite programming (SDP) problems. These are conic linear programs over the cone of positive semidefinite matrices,
formally defined as follows.

\begin{trailer}{Semidefinite programming (SDP) problem}
Assume we are given symmetric matrices $A_0,\ldots,A_m$ (all of the same size) and scalars $b_i\in\oR$ $(i \in [m])$. The  semidefinite programming problem in standard primal form is then defined
as
\[
p^*:=\inf_{X \succeq 0} \left\{\langle A_0,X\rangle \; : \; \langle A_i,X\rangle = b_i \; \forall i \in [m]\right\},
\]
where $\langle \cdot,\cdot\rangle$ now denotes the trace inner product, i.e., the Euclidean inner product  in the space of symmetric matrices,
and $X \succeq 0$ means that $X$ is a  symmetric  positive semidefinite matrix (corresponding to the L\"owner partial ordering of the symmetric matrices).

The dual semidefinite program reads
$$
d^*:= \sup_{y\in \oR^m} \left\{\sum_{i=1}^m b_iy_i: A_0-\sum_{i=1}^my_iA_i\succeq 0\right\}.
$$
Weak duality holds: $p^*\ge d^*$. Moreover, strong duality: $p^*=d^*$ holds, e.g., if the primal problem is bounded and admits a positive definite feasible solution $X$ (or if the dual is bounded and has a feasible solution $y$ for which $A_0-\sum_iy_iA_i$ is positive definite) (see, e.g., \cite{Barvinok_course in_convexity,BenTal}).

Next we recall how one can test whether a polynomial can be written as  a sum of squares of polynomials using semidefinite programming. This well known  fact plays a key role for reformulating the  inner and outer approximations of $\mathcal{M}(K)_+$ using semidefinite programs.
\end{trailer}

\newpage
\begin{trailer}{Checking sums of squares with SDP}

Given an integer $r\in \oN$ let 
$[x]_r=\{ x^\alpha: \alpha \in\oN^n_r\}$
consist of all monomials with degree at most $r$, thus the monomial basis of $\oR[x]_r$.

\begin{proposition}
\label{proposition:Gram matrix representation of sos}
For a given $n$-variate polynomial $h$, one has $h \in \Sigma_r$, if and only if the following polynomial identity holds:
\[
h(x) = [x]_r^\top M[x]_r \ \left(= \sum_{\alpha,\beta \in \mathbb{N}_r^n} M_{\alpha,\beta}x^{\alpha + \beta}\right),
\]
for some positive semidefinite matrix: $M = \left(M_{\alpha,\beta}\right)_{\alpha,\beta\in\oN^n_r} \succeq 0$. The above identity can be equivalently written as
\begin{equation}
\label{eq:h lmi}
h_\gamma = \sum_{\alpha,\beta \in \mathbb{N}_r^n \;:\; \alpha+\beta=\gamma} M_{\alpha,\beta} \quad \forall \gamma \in \mathbb{N}_{2r}^n.
\end{equation}
\end{proposition}
\end{trailer}

\paragraph{ SDP upper bounds for GPM via the inner  approximations}

Recall that the inner approximations of the cone $\mathcal{M}(K)_+$ restrict the measures on $K$  to the subsets $\mathcal{M}_{\mu_0}^r$ in \eqref{inner approximation}, i.e.\ to those
measures $\mu$ of the form $d\mu = h \cdot d\mu_0$, where $\mu_0$ is a fixed reference measure with $\text{Supp}(\mu_0)=K$ and $h \in \Sigma_r$ is a sum-of-squares polynomial density.

\ignore{
First note that the moments of such measure $\mu$  in $\mathcal{M}_{\mu_0}^r$ simplify to
\begin{eqnarray*}
  m^\mu_\alpha(K) &=& \int_K x^\alpha h(x)d\mu_0(x) \\
   &=& \int_K x^\alpha \sum_{\beta \in \mathbb{N}_{2r}^n} h_\beta x^\beta d\mu_0(x) \\
    &=& \sum_{\beta \in \mathbb{N}_{2r}^n} h_\beta \int_K   x^{\alpha+\beta} d\mu_0(x) \\
   &=& \sum_{\beta \in \mathbb{N}_{2r}^n} h_\beta m^{\mu_0}_{\alpha+\beta}(K).
\end{eqnarray*}
Thus the moments of $\mu$ are linear functions in the coefficients of $h$, involving the moments of $\mu_0$.
}

Replacing the cone $\mathcal{M}(K)_+$ in the GPM (\ref{prob:GPM}) by its subcone $\mathcal{M}^r_{\mu_0}$ we obtain the parameter
\begin{equation}\label{eq:lowr}
val_{inner}^{(r)} := \inf _{\mu\in \mathcal{M}^r_{\mu_0}} \left \{\int_K f_0(x)d\mu(x): \int_Kf_i(x)d\mu(x)=b_i\ \forall i\in [m]\right\},
\end{equation}
which provides  a hierarchy of upper bounds for GPM:
$$ val \le val_{inner}^{(r+1)}\le val_{inner}^{(r)}.
$$

According to the above discussion  these parameters can be reformulated as semidefinite programs involving the moments of the reference measure $\mu_0$. Indeed,
we may write the variable  density function as $h(x)=[x]_r^TM[x]_r$ with $M\succeq 0$ and arrive at the following  semidefinite program (in standard primal form).

\begin{trailer}{SDP formulation for the inner approximations based upper bounds }
\begin{equation}\label{eq:sdplowr}
val_{inner}^{(r)} =\inf_M \left\{ \langle A_0, M\rangle: \langle A_i, M\rangle =b_i \ \forall i\in [m], \ M=(M_{\alpha,\beta})_{\alpha,\beta\in\oN^n_r}\succeq 0\right\},
\end{equation}
where we set
$$
A_i=\int_Kf_i(x)[x]_r[x]_r^Td\mu_0(x)= \left(\int_K f_i(x)x^{\alpha +\beta}d\mu_0(x)\right)_{\alpha,\beta\in\oN^n_r} \ (0\le i\le m).
$$
Moreover, writing each polynomial $f_i$ in the monomial basis as
$f_i=\sum_{\gamma} f_{i,\gamma}x^\gamma$ one sees that the entries of the matrix $A_i$ depend linearly on the moments of the reference measure $\mu_0$, since
$\int_Kf_i(x)x^{\alpha+\beta}d\mu_0(x)= \sum_\gamma f_{i,\gamma} m^{\mu_0}_{\alpha+\beta+\gamma}(K).$
\end{trailer}

\ignore{
\begin{trailer}{Final inner approximation SDP formulation}
Thus the final SDP formulation of the inner approximation of the GPM \eqref{prob:GPM} looks as follows:
\begin{eqnarray*}
val_{inner}^{(r)} &:=& \min_{(h_\alpha)_{\alpha \in \mathbb{N}^n_{2r}}}  \sum_{\alpha \in \mathbb{N}_d^n} f_{0,\alpha}\sum_{\beta \in \mathbb{N}_{2r}^n} h_\beta m^{\mu_0}_{\alpha+\beta}(K) \\
&\mbox{s.t.}& \sum_{\alpha \in \mathbb{N}_d^n} f_{i,\alpha}\sum_{\beta \in \mathbb{N}_{2r}^n} h_\beta m^{\mu_0}_{\alpha+\beta}(K) = b_i \; \forall i \in [m], \\
&& h_\gamma = \sum_{\alpha,\beta \in \mathbb{N}_r^n \;:\; \alpha+\beta=\gamma} x_{\alpha,\beta} \quad \forall \gamma \in \mathbb{N}_{2r}^n\\
&&  X = \left(x_{\alpha,\beta}\right) \succeq 0.
\end{eqnarray*}
Note that one has $val_{inner}^{(r)} \ge val$ $(r = 1,2,\ldots)$.
\end{trailer}
}

To be able to compute the  above SDP one needs  the moments of the reference measure $\mu_0$ to be known on the set $K$. This is a restrictive assumption, since
even computing volumes of polytopes is an NP-hard problem. One is therefore restricted to specific choices of $\mu_0$ and $K$ where
the moments are known in closed form (or can be derived). In Table~\ref{tab:known moments} we therefore give an overview of some known moments for the Euclidean ball and sphere, the hypercube, and the standard simplex. (See \cite{Folland} for an easy derivation of the moments on the ball and the sphere.)
There we use the Gamma function: 
\begin{eqnarray*}
\Gamma(k)=(k-1)!,\quad  \Gamma\left(k+{1\over 2}\right)= \left(k-{1\over 2}\right) \left(k-1 -{1\over 2}\right) \cdots {1\over 2} \sqrt \pi \quad \text{ for  } k\in\oN.
\end{eqnarray*}

\ignore{
In the table we use the double factorial notation for an integer $k$:
\begin{eqnarray*}
k!! = \left\{ \begin{array}{l@{\quad}l}
k\cdot(k-2)\cdots 3\cdot1 & \textrm{if $k>0$ is odd,}\\
k\cdot(k-2)\cdots 4\cdot2 & \textrm{if $k>0$ is even,}\\
1 & \textrm{if $k=0$ or $k=-1$.}
\end{array} \right.
\end{eqnarray*}
}

\begin{table}[h!]
  \centering
\begin{tabular}{|c|c|}
   \hline
   $K$  & $m^{\mu_0}_{\alpha}(K)$  \\
    \hline
   $[0,1]^n$  & $\prod_{i=1}^n \frac{1}{\alpha_i+1}$ \\
   \hline
   $\Delta_n$  &  ${\prod_{i=1}^n\alpha_i!\over (\sum_{i=1}^n \alpha_i +n)!}$ \\
    \hline
    $S_n$ &
    $\left\{ \begin{array}{ll}
    {{2\Gamma(\beta_1)\cdots \Gamma(\beta_n)\over \Gamma(\beta_1+\ldots + \beta_n)}  }
    &    \mbox{if $\alpha \in (2 \mathbb{N})^n$} \quad \text{ with }\ \beta_i={\alpha_i+1\over 2}\ \text{ for } i\in [n]\\
0 & \mbox{otherwise} \\
\end{array} \right.$ \\
   \hline
   $B_n$ &
     $\left\{ \begin{array}{ll}
{1\over \alpha_1+\ldots +\alpha_n+n} {{2\Gamma(\beta_1)\cdots \Gamma(\beta_n)\over \Gamma(\beta_1+\ldots + \beta_n)}  }
&    \mbox{if $\alpha \in (2 \mathbb{N})^n$}\quad \text{ with } \ \beta_i={\alpha_i+1\over 2}\ \text{ for } i\in [n]\\
0 & \mbox{otherwise} \\
\end{array} \right.$ \\
   \hline
 \end{tabular}

  \caption{Examples of known moments for some choices of $K \subseteq \mathbb{R}^n$:
     $\Delta_n=\{x\in\oR_+^n: \sum_{i=1}^n x_i=1\}$ is the standard simplex and
   $B_n=\{x\in\oR^n:\|x\|\le 1\}$ is the unit Euclidean ball, in which case $\mu_0$ is the Lebesgue measure, and
   $S_n=\{x\in\oR^n: \|x\|=1\}$ is the unit Euclidean sphere in which case $\mu_0$ is the (Haar) surface measure on $S_n$.   }
  \label{tab:known moments}
\end{table}
If $K$ is an ellipsoid, one may obtain the moments of the Lebesgue measure on $K$ from the moments on the ball by an affine transformation of variables.
Also, if $K$ is a polytope, one may obtain the moments of the Lebesgue measure through triangulation of $K$, and subsequently using the formula for the simplex.

\paragraph{SDP lower bounds for GPM via the outer approximations}

Here we assume that $K$ is basic closed semi-algebraic, of the form
$$K=\{x\in\oR^n: g_j(x)\ge 0\ \ \forall j\in[k]\}, \quad \text{ where } g_1,\ldots,g_k\in\oR[x].$$
Recall that the dual cone of the truncated quadratic module generated by the polynomials $g_j$ describing the set  $K$  provides an  outer approximation of $ \mathcal{M}(K)_+$; we repeat its definition (\ref{dual:Qr})  for convenience:
 \[
\left(\mathcal{Q}^r(g_1,\ldots,g_k)\right)^* =  \left\{ \mu \in \mathcal{M}(K) \; : \; \int_K f d\mu \ge 0 \quad \forall f \in \mathcal{Q}^r(g_1,\ldots,g_k)\right\},
  \]
  where the quadratic module
 $\mathcal{Q}^r(g_1,\ldots,g_k)$ was defined in \eqref{eq:quadratic module}.

Replacing the cone $\mathcal{M}(K)_+$ in the GPM (\ref{prob:GPM}) by the above  outer approximations
we obtain the following parameters
\begin{equation}\label{eq:upr}
val_{outer}^{(r)} := \inf_{\mu \in \left(\mathcal{Q}^r(g_1,\ldots,g_k)\right)^*}
\left \{\int_K f_0(x)d\mu(x): \int_Kf_i(x)d\mu(x)=b_i\ \forall i\in [m]\right\},
\end{equation}
which provide a hierarchy of lower bounds for the GPM:
$$
val^{(r)}_{outer} \le val_{outer}^{(r+1)} \le val.
$$

Here too these parameters can be reformulated as semidefinite programs.
Indeed
a signed measure $\mu$ lies in the cone $\left(\mathcal{Q}^r(g_1,\ldots,g_k)\right)^*$ precisely when it  satisfies the condition
\begin{equation}\label{eq:cond}
\int_K g_j(x) \sigma_j(x)d\mu(x) \ge 0 \quad \forall \; \sigma_j \in \Sigma_{r_j}, \quad \forall j \in \{0,\ldots,k\},
\end{equation}
where $r_j=r-\lceil \deg(g_j)/2\rceil$.
Using Proposition \ref{proposition:Gram matrix representation of sos}, we may represent each
sum-of-squares $\sigma_j$ as 
\[
\sigma_j(x) = [x]_{r_j}^\top M^{(j)}[x]_{r_j} 
\]
for some matrix    $M^{(j)}\succeq 0$ (indexed by $\oN^n_{r_j}$). Hence we have
$$
\int_K g_j(x)\sigma_j(x)d\mu(x)= \int_K g_j(x)[x]_{r_j}^TM^{(j)}[x]_{r_j}d\mu(x) = \langle B_j^\mu,M^{(j)}\rangle,
$$ after setting
$$
B_j^\mu=\int_K g_j(x)[x]_{r_j}[x]_{r_j}^Td\mu(x) =\left(\int_K g_j(x)x^{\alpha+\beta}d\mu(x)\right)_{\alpha,\beta\in\oN^n_{r_j}}.
$$
Hence the condition (\ref{eq:cond}) can be rewritten as requiring, for each  $j\in\{0,1,\ldots,k\}$,
$$
\langle B_j^\mu,M^{(j)}\rangle \ge 0 \quad \text{ for all postive semidefinite matrices }\  M^{(j)} \text{ indexed by } \oN^n_{r_j},
$$
which in turn is equivalent to $B_j^\mu\succeq 0$
(since the cone of positive semidefinite matrices is self-dual).  Summarizing, the condition (\ref{eq:cond}) on the variable measure $\mu$ can be rewritten as
\begin{equation*}\label{eq:Bj}
B^\mu_j=\left(\int_K g_j(x)x^{\alpha+\beta}d\mu(x)\right)_{\alpha,\beta\in\oN^n_{r_j}}\succeq 0 \quad \forall j\in\{0,1,\ldots,k\}.
\end{equation*}
Finally, observe that only the moments of $\mu$ are playing a role in the above constraints. Therefore we may introduce new variables for these moments, say
$$
y_\alpha=\int_K x^\alpha d\mu(x)\quad \forall  \alpha\in\oN^n_{2r}.
$$
Writing the polynomials $g_j$ in the monomial basis as $g_j(x)=\sum_\gamma g_{j,\gamma}x^\gamma$ we arrive at the following SDP reformulation for the parameter
$val_{outer}^{(r)}.$

\begin{trailer} {SDP formulation for the outer approximations based lower bounds}
With $r_j=r-\lceil \deg(g_j)/2\rceil$ for $j\in\{0,1,\ldots,k\}$ and $d$ an upper bound on the degrees of $f_i$ for $i\in \{0,1,\ldots,m\}$ we have
\begin{eqnarray}
val_{outer}^{(r)}= \inf_{(y_\alpha)_{\alpha\in\oN^n_{2r}}} & \Big\{  \sum_{\alpha\in \oN^n_d}f_{0,\alpha} y_{\alpha}:  \sum_{\alpha\in \oN^n_d} f_{i,\alpha} y_{\alpha}=b_i\ \ \ \forall i\in [m],\label{eq:upr1} \\
&  \left(\int_K \sum_{\gamma} g_{j,\gamma} y_{\alpha+\beta+\gamma}     \right)_{\alpha,\beta\in\oN^n_{r_j}}\succeq 0 \quad \forall j\in\{0,1,\ldots,k\}\Big\}. \label{eq:upr2}
\end{eqnarray}
\end{trailer}

If, in the definition (\ref{eq:upr}) of  $val^{(r)}_{outer}$,  instead of the truncated quadratic module $\mathcal{Q}^r(g_1,\ldots,g_k)$ we use the larger quadratic module
$\mathcal{Q}^r(\prod_{j\in J}g_j: J\subseteq [k])$ generated by the pairwise products of the $g_j$'s, then we obtain a stronger bound on $val$, which we denote by
$\overline{val}^{(r)}_{outer}$. Thus
\begin{equation}\label{eq:upbar}
\overline{val}^{(r)}_{outer}
=\inf_{\mu \in (\mathcal{Q}^r(\prod_{j\in J}g_j: J\subseteq [k]))^*} \left\{\int_Kf_0(x)d\mu(x): \int_K f_i(x)d\mu(x) =b_i \ (i\in [m])\right\}
\end{equation}
and clearly we have
$$val^{(r)}_{outer} \le \overline{val}^{(r)}_{outer}\le val.$$
The parameter $\overline{val}^{(r)}_{outer}$ can also be reformulated as a semidefinite program, analogous to the program (\ref{eq:upr1})-(\ref{eq:upr2}), which however now involves $2^k+1$ semidefinite constraints instead of $k+1$ such constraints in (\ref{eq:upr2}) and thus its practical implementation is feasible only for small values of $k$.
On the other hand, as we will see later in Section \ref{sec:convPutinar}, the bounds $\overline{val}^{(r)}_{outer}$ admit a much sharper error analysis than the bounds $val^{(r)}_{outer}$ for the case of polynomial optimization.

\section{Convergence results for the inner approximation hierarchy}

In the rest of the paper we are interested in the convergence of the respective lower and upper SDP bounds on the optimal value of the GPM, as introduced in the previous
section. We will first consider in this section the upper bounds for the GPM arising from the inner approximations, since much more is known about their rate of convergence than for the lower bounds arising from the outer approximations. We deal first with the special case of polynomial optimization and then indicate how some of the results extend to the general GPM.

\subsection{The special case of global polynomial optimization}

Here we consider a special case of the GPM, namely global optimization of polynomials on compact sets (i.e.,\ problem \eqref{eq:global opt pol}) and review the main known results about the error analysis of the upper bounds $val_{inner}^{(r)}$.  After that in the next section  we will
explain how to extend this error analysis to the bounds for the general GPM problem.

Thus we now consider the problem
\begin{equation}\label{eq:minp}
val = \min_{x \in K} p(x),
\end{equation}
 asking to find the minimum value of the polynomial
$p(x) = \sum_{\alpha \in \mathbb{N}^n_d} p_\alpha x^\alpha$ over a compact set $K$.

Recall the definition of the inner approximation based upper bound (\ref{eq:lowr}), which can be rewritten here as
\[
val^{(r)}_{inner} = \min_{h \in \Sigma_r} \left\{ \int_K p(x)h(x) d\mu_0(x) \; : \; \int_K h(x) d\mu_0(x) = 1\right\},
\]
and its SDP reformulation from (\ref{eq:sdplowr}), which now reads
\begin{equation}\label{eq:lowrp}
val_{inner}^{(r)}=\min \left\{\langle A_0,M\rangle: \langle A_1,M\rangle=1,\ M=(M_{\alpha,\beta})_{\alpha,\beta\in\oN^n_r}\succeq 0\right\},
\end{equation}
with
$$A_0=\left(\int_Kp(x)x^{\alpha+\beta}d\mu_0(x)\right)_{\alpha,\beta\in\oN^n_r},\quad A_1=\left(\int_Kx^{\alpha+\beta}d\mu_0(x)\right)_{\alpha,\beta\in\oN^n_r},$$
where as before $\mu_0$ is a fixed reference measure on $K$.

A first observation made in \cite{Las11} is that this semidefinite program (\ref{eq:lowrp}) can in fact be reformulated as a generalized eigenvalue problem.
Indeed, its dual semidefinite program  reads
$$
\max\{\lambda: A_0-\lambda A_1\succeq 0\},$$
whose optimal value gives again the parameter $val_{inner}^{(r)}$ (since
strong duality holds).
Hence
 $val_{inner}^{(r)}$ is equal to
 the smallest generalized eigenvalue of the system
\begin{equation}
\label{eq:gen eig}
A_0v = \lambda A_1v, \quad v\ne 0.
\end{equation}
Thus one may compute $val_{inner}^{(r)}$ without having to solve an SDP problem.

In fact, if instead of the monomial basis $\{x^\alpha:\alpha \in \oR^n_{2r}\}$  we use
a polynomial basis $\{b_\alpha(x):\alpha\in\oN^n_{2r}\}$ of $\oR[x]_{2r}$ that is orthonormal with respect to the reference measure $\mu_0$
(i.e.,  such that $\int_K b_\alpha b_\beta d\mu_0 = 1 $ if $\alpha=\beta$ and 0 otherwise),  then
in the above semidefinite program (\ref{eq:lowrp}) we may set $A_1=I$ to be the identity matrix and
\begin{equation}\label{eq:A0}
A_0=\left(\int_Kp(x)b_\alpha(x)b_\beta(x)d\mu_0(x)\right)_{\alpha,\beta\in\oN^n_{2r}},
\end{equation}
whose entries now involve the `generalized' moments $\int_Kb_\alpha(x)d\mu_0(x)$ of $\mu_0$.
Then the parameter $val_{inner}^{(r)}$ can be computed as   the smallest eigenvalue of the matrix $A_0$:
\begin{equation}\label{eq:root}
val_{inner}^{(r)}= \lambda_{\min}(A_0)\quad \text{ where } A_0\   \text{ is as in (\ref{eq:A0})}.
\end{equation}

This fact was observed in \cite{KL_2018} and used there to establish a link with the roots of the orthonormal polynomials, permitting to analyze the quality of the bounds $val_{inner}^{(r)}$ for the case of the hypercube $K=[-1,1]^n$, see below for  details.


\medskip
In Table \ref{tab:convergence rates} we list the known convergence rates of the 
parameters $val_{inner}^{(r)}$ to the optimal value $val$ of problem \eqref{eq:minp},
i.e.,\ we review the known upper bounds for the sequence
$\{val_{inner}^{(r)} - val\}$, $r = 1,2,\ldots$

\begin{table}[h!]
\begin{center}
\begin{tabular}{|c|c|c|c|}
  \hline
  $K \subseteq \mathbb{R}^n$  & $val_{inner}^{(r)} - val$ & measure $\mu_0$ & reference\\ \hline
  compact & $o(1)$ & positive finite Borel measure & \cite{Las11}\\
  compact, satisfies interior cone condition &  $O\left(\frac{1}{\sqrt{r}}\right)$ & Lebesgue measure & \cite{KLS_MPA} \\
  convex body & $O\left(\frac{1}{r}\right)$ & Lebesgue measure & \cite{KL_MOR_2017}\\
 hypercube  $[-1,1]^n$ & $\Theta\left(\frac{1}{r^2}\right)$ & $ \prod_{i=1}^n (1-x_i^2)^{-1/2}dx_i$  &\cite{KL_2018}\\
  unit sphere, $p$ homogeneous & $O\left({1\over r}\right)$ & surface measure & \cite{DW}\\
     \hline
\end{tabular}
\caption{Known rates of convergence for the Lasserre hierarchy\label{tab:convergence rates} \ of upper bounds on $val$ in (\ref{eq:minp}) based on  inner approximations.}
\end{center}
\end{table}

We will give some details on the proofs of each of the four results listed in Table~\ref{tab:convergence rates}. After that we will mention an interesting connection with approximations based on cubature rules.


\paragraph{Asymptotic convergence}

   The first result in Table \ref{tab:convergence rates}  states that $\lim_{r \rightarrow \infty} val_{inner}^{(r)} = val$ if $K$ is compact and
  $\mu_0 \in \mathcal{M}(K)_+$. It is a direct consequence of the following result.

  \begin{theorem}[Lasserre \cite{Las11}]
Let ${K}\subseteq \mathbb{R}^n$ be compact, let $\mu_0$ be a fixed, finite, positive Borel measure 
with $\text{Supp}(\mu_0)={K}$.
 and let $f$ be a continuous function on $\mathbb{R}^n$. Then, $f$ is nonnegative on ${K}$ if and only if
\begin{equation*}
\int_{{K}}g^2fd\mu_0\ge0 \ \ \forall  g\in \oR[x].
\end{equation*}
\end{theorem}

The asymptotic convergence  of the bounds $val_{inner}^{(r)}$ to $val$ holds more generally for the minimization of a rational function $p(x)/q(x)$ over $K$ (assuming $q(x)>0$ for all $x\in K$).
Indeed, using the above theorem,  we obtain
\begin{eqnarray*}
\min_{x\in K}{p(x)\over q(x)}& =& \sup_{t\in \oR}\  t \ \text{ s.t. } p(x)\ge tq(x) \ \forall x\in K \\
&=&\sup_{t\in \oR}\ t \ \text{ s.t. }  \int_K p(x)h(x)d\mu_0(x)\ge t\int_{K} q(x)h(x)d\mu_0(x) \ \ \forall h\in \Sigma\\
& = & \inf_{h\in \Sigma} \int_K p(x)h(x)d\mu_0(x) \ \text{ s.t. } \int_Kq(x)h(x)d\mu_0(x)=1.
\end{eqnarray*}

\paragraph{Error analysis when $K$ is compact and satisfies an interior cone condition}

  The second result in Table \ref{tab:convergence rates} fixes the reference measure $\mu_0$ to the Lebesgue measure, and restricts the
  set $K$ to satisfy a so-called interior cone condition.
  \begin{definition}[Interior cone condition]
A set ${K}\subseteq \mathbb{R}^n$  satisfies an interior cone condition if there exist an angle $\theta\in (0,{\pi/ 2})$ and a radius $\rho>0$ such that,
 for every $x\in {K}$, a unit vector $\xi(x)$ exists such that
\[
\{x+\lambda y: y\in \mathbb{R}^n,\|y\|=1,y^T\xi(x)\ge\cos{\theta},\lambda\in[0,\rho] \} \subseteq {K}.
\]
\end{definition}
For example, all full-dimensional convex sets satisfy the interior cone condition for suitable parameters $\theta$ and $\rho$.
This assumption is used in \cite{KLS_MPA} to claim that the intersection of any ball with the set $K$ contains a positive fraction of the full ball, a fact used in the error analysis.

The main ingredient of the proof is  to
approximate the Dirac delta supported on a global minimizer by a Gaussian density of the form
\begin{equation}\label{eqH}
G(x) = \frac{1}{(2\pi \sigma^2)^{n/2}}\mbox{exp}\left(\frac{-\|x-x^*\|^2}{2\sigma^2}\right),
\end{equation}
where $x^*$ is a minimizer of $p$ on $K$,
and 
$\sigma^2 = \Theta(1/r)$.
Then we approximate the Gaussian density $G(x)$ by  a sum-of-squares polynomial $g_r(x)$ with degree
 $2r$. For this we use the  fact that the Taylor approximation of the exponential function $e^{-t}$
 is a sum of squares (since it is a univariate polynomial nonnegative on $\oR$).

 \begin{lemma}\label{lemexp}
For any $r\in \oN$ the univariate polynomial $\sum_{k=0}^{2r}{(-1)^k\over k!}  t^k$ (in the variable $t\in \oR$), defined as the Taylor expansion of the function $t\in \oR\mapsto e^{-t}$ truncated at degree $2r$,
 is a sum of squares of polynomials.
 \end{lemma}
Based on this the polynomial
 \[
g_r(x) = \frac{1}{(2\pi \sigma^2)^{n/2}} \sum_{k=0}^{2r}{(-1)^k \over k!} \left(\frac{-\|x-x^*\|^2}{2\sigma^2}\right)^k
\]
 is indeed a sum of squares with degree $2r$, which can be used (after scaling) as feasible solution within the definition of the bound $val_{inner}^{(r)}$.
 We refer to  \cite{KLS_MPA} for the details of the analysis.

\paragraph{Error analysis when $K$ is a convex body}

   The third item in Table \ref{tab:convergence rates} assumes that $K$ is now convex, compact and full-dimensional, i.e.,\ a convex body.
  The key idea is to use the following concentration result for the \emph{Boltzman density} (or \emph{Gibbs measure}).\footnote{This result
  is of independent interest in the study of \emph{simulated annealing} algorithms.}

  \begin{theorem}[Kalai-Vempala \cite{Kalai-Vempala 2006}]
If $p$ is a {linear} polynomial, $K$ is  a convex set, $T>0$ is a fixed `temperature' parameter, and $val = \min_{x \in K} p(x)$, then we have
\[
 \int_{\mathbf{K}}p(x)H(x)dx - val \leq  n T,
\]
where $$H(x) = \frac{\exp(-p(x)/T)}{\int_K \exp(-p(x)/T) dx}$$ is the Boltzman probability density supported on $K$.
\end{theorem}
The theorem still holds if $p$ is  convex, but not necessarily linear  \cite{KL_MOR_2017}. The proof of the
third item in Table \ref{tab:convergence rates} now proceeds as follows:

\begin{enumerate}
\item
Construct a sum-of-squares polynomial approximation $h_r(x)$ of the Boltzman density $H(x)$  by again using the fact that the even degree
{truncated Taylor expansion} of $e^{-t}$ 
 is a sum of squares (Lemma \ref{lemexp}); namely, consider the polynomial $h_r(x)= \sum_{k=0}^{2r}{(-1)^k\over k!} \left({-p(x)\over T}\right)^k$ (up to scaling).

\item
Use this construction to {bound the difference} between $val_{inner}^{(r)}$ and the Boltzman bound when choosing $T = O(1/r)$.

\item
Use the extension of the Kalai-Vempala result to get the required result for {convex} polynomials $p$.

\item
When  $p$ is nonconvex, the key ingredient is to reduce to the convex case by constructing a convex (quadratic)  polynomial $\hat p$ that upper bounds $p$ on $K$ and has the same minimizer on $K$, as indicated in the next lemma.

\begin{lemma}\label{lemquad}
Assume $x^*$ is a global minimizer of $p$ over $K$. Then the following  polynomial
$$\hat p(x)= p(x^*)+ \nabla p(x^*)^T(x-x^*)+ C_p \|x-x^*\|^2$$
with $C_p=\max_{x\in K} \|\nabla ^2p(x)\|_2$, is quadratic, convex, and separable. Moreover, it satisfies: $p(x)\le \hat p(x)$ for all $x\in K$, and $x^*$ is a global minimizer of $\hat p$ over $K$.
\end{lemma}

Then, in view of the inequality
\begin{equation}\label{eqkey}
\int_K \hat p h d\mu_0 \ge \int_K  p h d\mu_0 \quad  \forall h \in \Sigma_r,
\end{equation}
\end{enumerate}
it follows that the error analysis in the non-convex case follows directly from the error analysis in the convex case.
The details of the proof are given in \cite{KL_MOR_2017}.

\paragraph{Error analysis for the hypercube $K=[-1,1]^n$}

  The fourth result in Table \ref{tab:convergence rates} deals with the hypercube $K = [-1,1]^n$. A first key idea of the proof is that it suffices to show the $O(1/r^2)$ convergence
  rate for a univariate quadratic polynomial. This follows from Lemma \ref{lemquad} above (and (\ref{eqkey})), which implies that  it suffices to analyze the case of a quadratic, separable polynomial. Hence we may further restrict to the case when $K=[-1,1]$ and $p$ is a quadratic univariate polynomial.

\smallskip
In the univariate case, the key idea     is to use the eigenvalue reformulation of the bound $val_{inner}^{(r)}$ from (\ref{eq:root}). There, we use  the  polynomial basis $\{b_k: k\in\oN\}$ consisting of the  Chebyshev polynomials (of the first kind) which are orthonormal  with respect to the Chebyshev measure $d\mu_0$ on $K=[-1,1]$, indeed the measure used in  Table~\ref{tab:convergence rates}.

Then one may  use a connection to the extremal roots of these orthonormal polynomials.
Namely, for the linear polynomial $p(x)=x$, the parameter $val_{inner}^{(r)}$ coincides with the smallest root of the orthonormal polynomial $b_{r+1}$ (with degree $r+1$); this is a well known property of orthogonal polynomials, which follows from the fact that the matrix $A_0$ in (\ref{eq:A0}) is tri-diagonal and the 3-terms recurrence for the Chebyshev polynomials (see, e.g., \cite[\S1.3]{Dunkl}).
When $p$ is a quadratic polynomial, the matrix $A_0$ in the eigenvalue problem (\ref{eq:root})  is now 5-diagonal and `almost' Toepliz, properties that can be exploited to evaluate its  smallest eigenvalue.
 See
\cite{KL_2018} for details.

\paragraph{Error analysis for the unit sphere}

The last result in Table \ref{tab:convergence rates}  deals with the minimization of a homogeneous polynomial $p$ over the unit sphere 
$S_n=\{x\in\oR^n: \sum_{i=1}^n x_i^2=1\}$, in which case Doherty and Wehner \cite{DW} show a convergence rate in $O(1/r)$.
Their construction for a suitable sum-of-squares polynomial density in $\Sigma_r$ is in fact closely related to their analysis of the outer approximation based lower bounds $val^{(r)}_{outer}$. Doherty and Wehner \cite{DW} indeed show the following stronger result:
$val^{(r)}_{inner}-val^{(r)}_{outer} =O(1/r)$, to which we will come back  in Section \ref{sec:convPutinar} below.

\paragraph{Link with positive cubature rules }
There is an interesting link between positive cubature formulas and the upper bound
\[
val^{(r)}_{inner} = \min_{h \in \Sigma_r} \left\{ \int_K ph d\mu_0 \; : \; \int_K h d\mu_0 = 1\right\},
\]
which was recently pointed out in \cite{Martinez et al quadrature} and is summarized in the next result.

\begin{theorem}[Martinez et al.\ \cite{Martinez et al quadrature}]
\label{thm:quadrature}
Let $x^{(1)},\ldots,x^{(N)} \in K$ and weights $w_1 >0,\ldots,w_N >0$ give a positive
cubature rule on $K$ for the measure $\mu_0$, that is exact for polynomials of
total degree at most $d + 2r$, where $d > 0$ and $r > 0$ are given
integers. Let $p$ be a polynomial of degree $d$.

Then, if $h$ is a polynomial nonnegative on $K$ and of degree at most $2r$, one has
\[
 \int_K ph d\mu_0 \ge \min_{\ell \in [N]} p(x^{(\ell)}).
\]

In particular, the inner approximation bounds therefore satisfy
 \[
val^{(r)}_{inner} \ge \min_{\ell \in [N]} p(x^{(\ell)}).
\]
\end{theorem}
The proof is an immediate consequence of the definitions, but this result has several interesting implications.

\begin{itemize}
\item
First of all, one may derive information about the rate of convergence  for the scheme
$\min_{\ell \in [N]} p(x^{(\ell)})$ from  the error bounds in Table \ref{tab:convergence rates}.
For example, if $K$ is a convex body, the implication is that $\min_{\ell \in [N]} p(x^{(\ell)}) - val= O(1/r)$.

\smallskip
\item
Also, if a positive cubature rule is known for the pair $(K,\mu_0)$, and the number of points $N$ meets the Tchakaloff bound
$N = {n+2r+d \choose 2r+d}$, then there is no point in computing the parameter $val^{(r)}_{inner}$. Indeed,  as $$val^{(r)}_{inner} \ge \min_{\ell \in [N]} p(x^{(\ell)})\ge val,$$  the  right-hand-side bound is stronger and can  be computed more efficiently.
Having said that, positive cubature rules that meet the Tchakaloff bound are only known in special cases, typically in low dimension and degree; see e.g.\ \cite{Cools Encyclopedia cubature,Trefethen cubature review,Dai-Xu book}, and the references therein.

\smallskip\item
Theorem \ref{thm:quadrature} also shows why the last convergence rate in Table \ref{tab:convergence rates} is tight for $K = [-1,1]^n$.
Indeed if we consider the univariate example $p(x) = x$ and the Chebyshev probability measure $d\mu_0(x) = \frac{1}{\pi\sqrt{1-x^2}}dx$ on $K=[-1,1]$, then a positive cubature scheme is given
by
\[
x^{(\ell)} = \cos \left(\frac{2\ell-1}{2N}\pi\right), \; w_\ell= \frac{1}{N} \quad \forall \ell \in [N],
\]
and it is exact at degree $2N-1$.
This is known as the Chebyshev-Gauss quadrature, and the points are precisely the roots of the degree $N$ Chebyshev polynomial of the first kind.
Thus, with $N=r+1$, in this case we have
 \[
val^{(r)}_{inner} \ge \min_{\ell \in [N]} p(x^{(\ell)}) = \min_{\ell\in [N]} \cos \left(\frac{(2\ell-1)\pi}{2N}\right) = \cos \left( - \pi/(2N)\right) =-1 + \Omega\left(\frac{1}{N^2}\right).
\]
This explains that the $\Theta(1/r^2)$ result in Table \ref{tab:convergence rates}  holds for $p(x)=x$.
A different proof of this result is given in \cite{KL_2018}, where it is shown that for this example one actually has equality $val^{(r)}_{inner} =  \cos \left( - \pi/(2N)\right)$.

\smallskip\item
Finally, Theorem \ref{thm:quadrature} shows that there is not much gain
in using a larger set of densities than $\Sigma_r$  in the definition of the inner approximations $\mathcal{M}^r_{\mu_0}$
since the statement of the theorem
holds for any nonnegative polynomial $h$ on $K$. For example, for the hypercube $K=[-1,1]^n$, if we use  the larger set of densities $h\in\mathcal{Q}^r(\prod_{j\in J}(1-x_j^2): J\subseteq [k])$ and the Chebyshev measure as reference measure $\mu_0$ on $[-1,1]^n$,  then we obtain upper bounds with  convergence rate in $O(1/r^2)$  \cite{DHL SIOPT}.
This  also follows from the later results in \cite{KL_2018} where in addition  it is shown that
this convergence result is tight for linear polynomials. By the above discussion tightness also follows from Theorem \ref{thm:quadrature}.
\end{itemize}

\paragraph{{Upper bounds using grid point sets}}

Of course one may also obtain upper bounds on $val$, the minimum value taken by a polynomial $p$ over a compact set $K$, by evaluating $p$ at any suitably selected set of points in $K$. This corresponds to restricting the optimization over selected   finite atomic measures in the definition of $val$.

A first basic idea is to select the grid point sets consisting of all rational points in $K$ with denominator $r$ for increasing values of $r\in\oN$.
For the standard simplex $K=\Delta_n$ and the hypercube $K=[0,1]^n$ this leads to upper bounds that satisfy:
\begin{equation}\label{eq:simplex1}
\min_{x\in K, rx\in\oN^n}p(x)- \min_{x\in K}p(x)\le {C_d\over r}\left(\max_{x\in K}p(x)-\min_{x\in K}p(x)\right) \quad \text{ for all } { r\ge d,}\end{equation}
where $C_d$ is a constant that depends only on the degree $d$ of $p$;
see \cite{KLP} for $K=\Delta_n$ and  \cite{KL10} for $K=[0,1]^n$.
A faster regime in $O(1/r^2)$ can be shown when allowing a constant that depends on the polynomial $p$ (see \cite{KLSV} for $\Delta_n$ and \cite{KLLS17} for $[0,1]^n$).
Note that  the number of  rational points with denominator $r$ in the simplex $\Delta_n$ is  ${n+r-1\choose r}=O(n^r)$  and thus the computation time for these upper bounds is polynomial in the dimension $n$ for any fixed order $r$. On the other hand,
there are   $(r+1)^n=O(r^n)$ such grid points in  the hypercube $[0,1]^n$ and thus the computation time of the upper bounds grows exponentially with the dimension $n$.

\ignore{
 instance, when $K=\Delta_n$ is the standard simplex
we may select the grid point sets
consisting of the rational points in $\Delta_n$ with denominator $r$ for increasing values of $r\in \oN$.
For the simplex and the hypercube the following error bounds hold:
\begin{equation}\label{eq:simplex1}
\min_{x\in\Delta_{n}, rx\in\oN^n}p(x)- \min_{x\in \Delta_n}p(x)\le {C_d\over r}\left(\max_{x\in\Delta_n}p(x)-\min_{x\in \Delta_n}p(x)\right)
\end{equation}
where the constant $C_d$ depends only on the degree $d$ of $p$ \cite{KLP}. As the number of such rational points is  ${n+r-1\choose r}=O(n^r)$  the computation time for the upper bounds is polynomial in the dimension $n$ for any fixed order $r$.
Note that one can show an faster regime in $O(1/r^2)$ when allowing a constant depending on $p$:
$$\min_{x\in\Delta_{n}, rx\in\oN^n}p(x)- \min_{x\in \Delta_n}p(x)\le {C_p\over r^2},$$
where the constant $C_p$ now depends on $p$ \cite{KLSV}.

For the hypercube $K=[0,1]^n$ if we select all grid points with denominator $r$ in $[0,1]^n$ then the same regime in $O(1/r^2)$ holds:
$$\min_{x\in [0,1]^n, rx\in \oN^n} p(x) -\min_{x\in [0,1]^n} p(x)\le {C_p\over r^2}$$
for some constant $C_p$ depending on $p$ \cite{KLLS17}. However the computation of this upper bound needs $O(r^n)$ evaluations which is now exponential in $n$ for any fixed order $r$.
}

For a general convex body $K$ some constructions are proposed recently in \cite{PV} for suitable grid point sets (so-called meshed norming sets)  $X_d(\epsilon)\subseteq K$ where $d\in\oN$ and $\epsilon>0$. Namely, whenever $p$ has degree at most $d$, by minimizing $p$ over $X_d(\epsilon)$ one obtains an upper bound on the minimum of $p$ over $K$ satisfying
$$\min_{x\in X_d(\epsilon)}p(x)-\min_{x\in K}p(x) \le \epsilon\left(\max_{x\in K}p(x)-\min_{x\in K}p(x)\right),$$
where the computation involves $|X_d(\epsilon)|= O\left(\left({d\over \sqrt \epsilon}\right)^{2n}\right)$ point evaluations, thus  exponential in the dimension $n$ for fixed precison $\epsilon$.

In comparison, the computation of the upper bound $val^{(r)}_{outer}$ relies on a semidefinite program involving a matrix of size ${n+r\choose r}=O(n^r)$, which is polynomial in the dimension $n$ for any fixed order $r$.

\subsection{The general problem of moments (GPM)}
One may extend the results of the last section to the inner approximations for the general GPM \eqref{prob:GPM}.
In other words, we now consider the upper bounds (\ref{eq:lowr})  obtained using the inner approximations of the cone $\mathcal{M}(K)_+$, which we repeat for convenience:
\[
val^{(r)}_{inner} = \inf_{h \in \Sigma_r} \left\{ \int_K f_0(x)h(x) d\mu_0(x) \; : \; \int_K f_i(x) h(x)d\mu_0(x) = b_i \quad \forall i \in [m]\right\}.
\]
A first observation is that this program 
may not have a feasible solution, even if the GPM \eqref{prob:GPM} does.
For example,  two constraints like
\[
\int_0^1 x d\mu(x) = 0,\ \int_0^1 d\mu(x)=1
\]
admit the Dirac measure $\mu=\delta_{\{0\}}$ as solution but they do not admit any solution of the form
$d\mu = h dx$ with $h \in \Sigma_r$ for any $r \in \mathbb{N}$.
Thus any convergence result must relax the equality constraints of the GPM \eqref{prob:GPM} in some way, or involve additional assumptions.

We now  indicate how one may use the convergence results of the last section to derive an error analysis for  the  inner approximations of the GPM when relaxing the equality constraints.

\begin{theorem}[De Klerk-Postek-Kuhn \cite{DK-P-K dist robust opt}]
\label{th:moments sos}
Assume that $f_0,\ldots,f_m$ are polynomials, $K$ is compact and  the GPM \eqref{prob:GPM} has an optimal solution.
Let $b_0:=val $ denote the optimum value of (\ref{prob:GPM}) and for
any integer $r\in \oN$ define the parameter
$$\Delta^{(r)}:= \min_{h\in\Sigma_r} \max_{ i\in \{0,1,\ldots,m\}}\left| \int_K f_i(x)h(x)d\mu_0(x)-b_i\right|.$$
Then the following assertions hold:
\begin{enumerate}
\item[(1)] $\displaystyle\lim_{r\to \infty} \Delta^{(r)} = 0$.

\item[(2)]
$\Delta^{(r)}=O\left({1\over r^{1/4}}\right)$
 if $K$ satisfies an interior cone assumption and $\mu_0$ is the Lebesgue measure;
\item[(3)]
$\Delta^{(r)}=O\left({1\over r^{1/2}}\right)$
if $K$ is a convex body and $\mu_0$ is the Lebesgue measure;
\item[(4)] $\Delta^{(r)}=O\left({1\over r}\right)$
if $K = [-1,1]^n$  and $d\mu_0(x) = \prod_i (1-x_i^2)^{-1/2}dx_i$.
\end{enumerate}
\end{theorem}

We will derive this from the convergence results for global polynomial optimization in Table \ref{tab:convergence rates}. By assumption,  problem (\ref{prob:GPM}) has an optimal solution and by Theorem~\ref{thm:Rogosinsky} we may assume it has an atomic optimal solution $\mu^*=\sum_\ell\lambda_\ell\delta_{x^*_\ell}$ with $\lambda_\ell>0$ and $x^*_\ell\in K$.
We now sketch the proof.

\begin{enumerate}
  \item For each atom  $x^*_\ell$  of the optimal measure $\mu^*$  consider  the  polynomial
  \[
  p_\ell(x) = \sum_{i=0}^m \left(f_i(x) - f_i(x^*_\ell)\right)^2,
  \]
  whose minimum value over $K$ is equal to 0 (attained at $x^*_\ell$).
  \item
 We  apply the error analysis
  of the previous section to the problem of minimizing  the polynomial $p_{\ell}$ over $K$. In particular, the asymptotic convergence of the upper bounds implies that for any given $\epsilon>0$
  $$
 \exists r\in\oN \ \ \exists h_{\ell}\in\Sigma_{r}\  \text{ s.t.} \int_K p_{\ell}(x) h_{\ell}(x)d\mu_0(x) \le \epsilon^2, \int_K h_{\ell}(x)d\mu_0(x)=1
  $$
   and, therefore,
    \begin{equation}\label{eqfi}
     \int_K (f_i(x) -f_i(x^*_\ell))^2h_{\ell}(x)d\mu_0(x)  \le \epsilon^2\ \ \ \forall  i \in \{0,\ldots,m\}.
     \end{equation}
  \item
Using the Jensen inequality, one obtains
\[
\left|\int_K f_i(x)h_{\ell}(x)d\mu_0(x) - f_i(x^*_\ell)\right| = \left|\int_K (f_i(x) -f_i(x^*_\ell))h_{\ell}(x)d\mu_0(x)\right|   \le \epsilon 
\]
for each $i\in \{0,\ldots,m\}$.
\item
We now consider
the sum-of-squares density $h := \sum_\ell \lambda_\ell h_{\ell} \in \Sigma_r$. Then we have $b_i=\int_K f_i(x)d\mu^*(x)=\sum_\ell \lambda_\ell f_i(x^*_\ell)$ for each $i\in\{0,\ldots,m\}$.
Moreover,
 the above argument shows that for any $i\in\{0,\ldots,m\}$
\[
 \left|\int_K f_i(x)h(x)d\mu_0(x) - b_i \right| =
 \left|\sum_\ell\lambda_\ell\left(\int_K f_i(x)h_{\ell}(x)d\mu_0(x) - f_i(x^*_\ell)\right) \right|
 \le \epsilon \mu^*(K)
\]
with $\mu^*(K)=\sum_\ell\lambda_\ell$.
This shows that $\Delta^{(r)}\le \epsilon\mu^*(K)$ and thus  the desired asymptotic result (1).
\item
The additional three claims (2)-(4) follow in the same way using the results in Table \ref{tab:convergence rates}. For instance, in case (1) when $K$ satisfies an interior cone condition and $\mu_0$ is the Lebesgue measure, we replace the estimate (\ref{eqfi}) by
$$\left|\int_K(f_i(x)-f_i(x^*_\ell))^2 h_\ell(x)d\mu_0(x)\right|=O\left({1\over \sqrt r}\right),$$ which leads to  $\Delta^{(r)} = O\left({1\over r^{1/4}}\right)$ (since we `loose a square root' when applying Jensen inequality).
\end{enumerate}

\medskip
We may also use the relation  with positive cubature rules discussed in the previous section (Theorem \ref{thm:quadrature}) to obtain the following cubature-based approximations for the GPM (\ref{prob:GPM}).

\ignore{
\begin{corollary}
Assume a positive cubature rule for $(K,\mu_0)$ that is exact for polynomials of degree at most $d+2r$ is given
by points $x^{(i)} \in K$ with corresponding weights $w_i \ge 0$ for $i \in [N]$.
Under the assumptions of Theorem \ref{th:moments sos}, there exists  $\lambda_i \ge 0$ ($i \in [N]$) such that, for $\mu^* = \sum_{i=1}^N \lambda_i \delta_{x^{(i)}}$,
one has
\[
\left| \int_K f_0(x) d\mu^*(x) - val \right| \le \epsilon_r, \; \left| \int_K f_i(x)  d\mu^*(x) - b_i \right| \le \epsilon_r \quad \forall i \in [m],
\]
where $\epsilon_r > 0$ satisfies the following conditions:
 \begin{enumerate}
 \item
 $\lim_{r \rightarrow \infty} \epsilon_r = 0$;
\item
$\epsilon_r =   O(r^{-1/4})$ if $K$ satisfies interior cone assumption and $\mu_0$ is the Lebesgue measure;
\item
$\epsilon_r =   O(r^{-1/2})$ if $K$ is a convex body and $\mu_0$ is the Lebesgue measure;
\item
$\epsilon_r =   O(r{-1})$ if $K = [-1,1]^n$  and $d\mu_0(x) = \prod_i (1-x_i^2)^{-1/2}dx_i$;
\end{enumerate}
\end{corollary}
Note that the $\lambda_i$ in the corollary may be found using linear programming, for given $r$ and $\epsilon_r$, if the positive cubature rule
is known. (A similar idea was used in \cite{Ryu_Boyd_quadrature}.)
}

\begin{corollary}\label{cor:cubatureGPM}
Assume the GPM (\ref{prob:GPM}) admits an optimal solution and let $d$ denote the maximum degree of the polynomials $f_0,\ldots,f_m$.
For any integer $r\in\oN$ assume we have a cubature rule for $(K,\mu_0)$ that is exact for degree $d+2r$, consisting of the points $x^{(\ell)}\in K$ and weights $w_\ell>0$ for $\ell\in [N]$, and define the parameter
$$ \Delta^{(r)}_{cub}:=\min_\nu \max_{ i\in \{0,1,\ldots,m\}} \left|\int_K f_i(x)d\nu -b_i\right|, $$
where in the outer minimization we minimize over all atomic measures $\nu$  whose atoms all belong to the set $ \{x^{(\ell)}: \ell\in [N]\}.$
Then the following assertions hold:

 \begin{enumerate}
 \item[(1)]
 $\displaystyle\lim_{r \rightarrow \infty}  \Delta^{(r)}_{cub} = 0$;
\item[(2)]
$ \Delta^{(r)}_{cub}=O\left({1\over r^{1/4}}\right)$  if $K$ satisfies an interior cone assumption and $\mu_0$ is the Lebesgue measure;
\item[(3)]
$ \Delta^{(r)}_{cub}=O\left({1\over \sqrt r}\right)$  if $K$ is a convex body and $\mu_0$ is the Lebesgue measure;
\item
$ \Delta^{(r)}_{cub}=O\left({1\over r}\right)$ if $K = [-1,1]^n$  and $d\mu_0(x) = \prod_i (1-x_i^2)^{-1/2}dx_i$.
\end{enumerate}
\end{corollary}
This result follows  from Theorem \ref{th:moments sos}. Indeed, for any polynomial $h\in\Sigma_r$, the polynomials $f_ih$ have degree at most $d+2r$ so that using the cubature rule we obtain
$$\int_Kf_i(x)h(x)d\mu_0(x) = \sum_{\ell=1}^N w_\ell f_i(x^{(\ell)}) h(x^{(\ell)})= \int_K f_i(x) d\nu(x),$$ where $\nu$ is the atomic measure with atoms $x^{(\ell)}$ and weights $\alpha_\ell:=w_\ell h(x^{(\ell)})$ for $\ell\in [N]$. Therefore, the parameter  $\Delta ^{(r)}_{cub}$ in Corollary \ref{cor:cubatureGPM} is upper bounded by the parameter
$ \Delta^{(r)}$ in Theorem \ref{th:moments sos}.
The claims (1)-(4) now follow directly from the corresponding claims in Theorem \ref{th:moments sos}.

Note that, for any fixed $r\in \oN$, in order to find the best atomic measure $\nu$ in the definition of $ \Delta^{(r)}_{cub}$ we need to find the best  weights $\alpha_\ell$ ($\ell\in [N]$)  giving the measure  $\nu=\sum_{\ell=1}^N \alpha_\ell\delta_{x^{(\ell)}}$. This can be done   by solving the following linear program:
$$\Delta^{(r)}_{cub}=\min_{t,\alpha_\ell \in\oR} \ t\  \text{ s.t. } \alpha_\ell\ge 0\ (\ell\in [N]),\ \left |\sum_{\ell=1}^N \alpha_\ell f_i(x^{(\ell)}) -b_i \right|\le t\ \  \forall i\in \{0,1,\ldots,m\}.$$
(This is similar to an idea used in  \cite{Ryu_Boyd_quadrature}.)

\section{Convergence results for the outer approximations}

In this last section  we consider the convergence of the lower bounds 
 for the GPM (\ref{prob:GPM}), that  are obtained by using  outer approximations for the cone of positive measures.
We first mention properties dealing with asymptotic and finite convergence for the general GPM and after that we mention some known results on the error analysis in the special case of polynomial optimization.

Here we assume $K$ is a compact semi-algebraic set, defined as before by
$$K=\{x\in\oR^n: g_j(x)\ge 0\quad \forall j\in [k]\},$$
where $g_1,\ldots,g_k\in \oR[x].$ 
We will consider the following  ({\em Archimedean}) condition:
\begin{equation}\label{eq:ass}
\exists r\in \oN \ \exists u\in \mathcal{Q}^r(g_1,\ldots,g_k)\ \text{ s.t.  the set } \{x\in \oR^n: u(x)\ge 0\} \text{ is compact.}
\end{equation}
This condition clearly implies that $K$ is compact. Moreover, it  does not depend  on the set $K$ but  on the choice of the polynomials used to describe $K$.
Note  that it is easy to modify the presentation of $K$ so that the condition (\ref{eq:ass}) holds. Indeed, if we know the radius $R$ of a ball containing $K$
 then, by adding to the description of $K$ the (redundant) polynomial constraint  $g_{k+1}(x):= R^2-\sum_{i=1}^n x_i^2 \ge 0$, we can ensure that assumption (\ref{eq:ass}) holds for this enriched presentation of $K$.

For convenience we recall the definition of the bounds $val^{(r)}_{outer}$ from (\ref{eq:upr}):
\begin{equation*}
val_{outer}^{(r)} = \inf_{\mu \in \left(\mathcal{Q}^r(g_1,\ldots,g_k)\right)^*}
\left \{\int_K f_0(x)d\mu(x): \int_Kf_i(x)d\mu(x)=b_i\ \forall i\in [m]\right\},
\end{equation*}
where  we refer to (\ref{eq:quadratic module}) and (\ref{dual:Qr}) for the definitions of the truncated quadratic module $\mathcal{Q}^r(g_1,\ldots,g_k)$ and of its dual cone  $(\mathcal{Q}^r(g_1,\ldots,g_k))^*$.

We also recall the stronger bounds $\overline{val}^{(r)}_{outer}$,  introduced in (\ref{eq:upbar}), and obtained by replacing in the definition
of $val_{outer}^{(r)}$ the  cone $\mathcal{Q}^r(g_1,\ldots,g_k)$
by the larger cone $\mathcal{Q}^r(\prod_{j\in J}g_j: J\subseteq [k]))$, so that we have
$$val_{outer}^{(r)}\le \overline{val}^{(r)}_{outer}\le val.$$

\ignore{
We will also briefly consider  the stronger bounds on $val$, denoted by $\overline{val}^{(r)}_{outer}$, that are obtained by replacing in the above definition of $val^{(r)}_{outer}$ the truncated quadratic module $\mathcal{Q}^r(g_1,\ldots,g_k)$ by the larger cone
$\mathcal{Q}^r\left(\prod_{j\in J}g_j: J\subseteq [k]\right)$,
thus the truncated quadratic module generated by all pairwise products of the $g_j$'s.
Clearly we have
$$val^{(r)}_{outer} \le \overline{val}^{(r)}_{outer}\le val.$$
This new parameter can also be reformulated as a semidefinite program, analogous to the program (\ref{eq:upr1})-(\ref{eq:upr2}), but it now involves $2^k+1$ semidefinite constraints instead of $k+1$ such constraints in (\ref{eq:upr2}) and thus its practical implementation is feasible only for small values of $k$.
On the other hand, as we will see later in this section the bounds $\overline{val}^{(r)}_{outer}$ admit a much sharper error analysis than the bounds $val^{(r)}_{outer}$ for the case of polynomial optimization.
}

\subsection{Asymptotic and finite convergence}
Here we present some results on the asymptotic and finite convergence of the lower bounds  on $val$ obtained by considering outer approximations of the cone $\mathcal{M}(K)_+$.

\paragraph{Asymptotic convergence}
The parameters $val^{(r)}_{outer}$ form a non-decreasing sequence of lower bounds for the optimal value $val$ of problem (\ref{prob:GPM}), which converge to it under
 assumption (\ref{eq:ass}). This asymptotic convergence result relies on the following representation result of Putinar \cite{Putinar} for positive polynomials.

 \begin{theorem}[Putinar]\label{thm:Putinar}
 Assume $K$ is compact and assumption (\ref{eq:ass}) holds.
Any polynomial $f$ that is strictly positive on $K$ (i.e., $f(x)>0 $ for all $x\in K$) belongs to  $\mathcal{Q}^r(g_1,\ldots,g_k)$ for some  $r\in\oN$.
 \end{theorem}

  The following result can be found in  \cite{Las2008, lass-book} for the general GPM and in \cite{lass-siopt-01} for the case of global polynomial optimization.

\begin{trailer}{Asymptotic convergence for the bounds $val^{(r)}_{outer}$}
  \begin{theorem} \label{thm:asympconv}
  Assume $K$ is compact and assumption (\ref{eq:ass}) holds. Then we have
  $$val^*\le lim_{r\to\infty} val^{(r)}_{outer} \le val,$$
  with equality: $val^*= lim_{r\to\infty} val^{(r)}_{outer} = val$
   if, in addition, there exists $z\in\oR^{m+1}$ such that $\sum_{i=0}^m z_if_i(x)>0$ for all $x\in K$.
  \end{theorem}
\end{trailer}

This result follows using Theorem \ref{thm:Putinar}.
Observe  that it  suffices to show the inequality: $val^*\le \sup_r val^{(r)}_{outer}$ (as the rest follows using Corollary \ref{thm:sduality2}).
For this let $\epsilon>0$ and let $y\in \oR^m$ be feasible for $val^*$, i.e., $f_0(x)-\sum_{i=1}^m y_if_i(x)\ge 0$ for all $x\in K$; we will show the inequality
$b^Ty\le \sup_r val^{(r)}_{outer} +\epsilon\mu(K).$ Then, letting $\epsilon$ tend to 0 gives $b^Ty \le \sup_r val^{(r)}_{outer}$ and thus the desired result: $val^*\le \sup_r val^{(r)}_{outer}=\lim_{r\to\infty}val^{(r)}_{outer}$.

As the polynomial $f_0+\epsilon -\sum_iy_if_i$ is strictly positive on $K$, it belongs to $\mathcal{Q}^r(g_1,\ldots,g_k)$ for some $r\in\oN$ in view of Theorem \ref{thm:Putinar}. 
Then, for any measure $\mu$ feasible for $val^{(r)}_{outer}$, we have
$\int_K(f_0+\epsilon-\sum_iy_if_i)d\mu\ge 0$, which implies $b^Ty\le \int_Kf_0d\mu + \epsilon \mu(K)$ and thus the desired inequality:
$$b^Ty\le  val^{(r)}_{outer} +\epsilon \mu(K) \le \sup_r val^{(r)}_{outer}+\epsilon \mu(K).$$

\medskip
When assuming only $K$ compact (thus not assuming condition (\ref{eq:ass})),  the following representation result of Schm\"udgen \cite{Schm} permits to show
the asymptotic convergence of the  stronger  bounds $\overline{val}^{(r)}_{outer}$ to $val$ (in the same way as  Theorem \ref{thm:asympconv} follows from Putinar's theorem).

\begin{theorem}[Schm\"udgen] \label{thm:Schmudgen}
 Assume $K$ is compact.
Any polynomial $f$ that is strictly positive on $K$ (i.e., $f(x)>0 $ for all $x\in K$) belongs to
$\mathcal{Q}^r(\prod_{j\in J}g_j: J\subseteq [k])$ for some $r\in\oN$.
\end{theorem}

\begin{trailer}{Asymptotic convergence for the bounds $\overline{val}^{(r)}_{outer}$}
\begin{theorem}
  Assume $K$ is compact. Then we have
  $$val^*\le lim_{r\to\infty} \overline{val}^{(r)}_{outer} \le val,$$
  with equality: ${val}^*= lim_{r\to\infty} \overline{val}^{(r)}_{outer} = val$
   if, in addition, there exists $z\in\oR^{m+1}$ such that  $\sum_{i=0}^m z_if_i(x)>0$ for all $x\in K$.
\end{theorem}
\end{trailer}

\paragraph{Finite convergence}

A remarkable property of the lower bounds $val^{(r)}_{outer}$ is that they often exhibit finite convergence.
Indeed, there is an easily checkable criterion, known as the {\em flatness condition}, that permits to conclude that the bound is exact: $val^{(r)}_{outer}=val$, and to extract  an (atomic) optimal solution to the GPM. This is condition (\ref{eq:flat}) below, which permits to claim that a given truncated sequence  is indeed the sequence of moments of a positive measure; it goes back to work of Curto and Fialkow (\cite{CF}, see also \cite{lass-book,Lau09} for details). To expose it we use the SDP formulation (\ref{eq:upr1})-(\ref{eq:upr2}) for the parameter $val^{(r)}_{outer}$.

\begin{trailer}{Finite convergence}
\begin{theorem} (see \cite[Theorem 4.1]{lass-book})
Set $d_K:=\max\{\lceil\deg(g_j/2\rceil: j\in [k]\}$ and let  $r\in \oN$ such that $2r\ge \max\{\deg(f_i): i\in \{0,\ldots,m\}\}$ and $r\ge d_K$.
Assume the program (\ref{eq:upr1})-(\ref{eq:upr2}) defining the parameter $val^{(r)}_{outer}$ has an optimal solution $y=(y_\alpha)_{\alpha\in\oN^n_{2r}}$ that satisfies the following (flatness) condition:
\begin{equation}\label{eq:flat}
\text{rank} M_s(y)=\text{rank} M_{s-d_K}(y) \quad \text{ for some integer } s \text{ s.t. } d_K\le s\le r,
\end{equation}
where
$$M_s(y)=(y_{\alpha+\beta})_{\alpha,\beta\in\oN^n_s} \ \text{ and } \ M_{s-d_K}(y)=(y_{\alpha+\beta})_{\alpha,\beta\in\oN^n_{s-d_K}}.
$$
Then  equality $val^{(r)}_{outer}=val$ holds and the GPM problem (\ref{prob:GPM}) has an optimal solution $\mu\in\mathcal{M}(K)_+$ which is atomic and supported on $\text{rank} M_s(y)$ points in $K$.
\end{theorem}
\end{trailer}

Under the flatness condition (\ref{eq:flat}) there is an algorithmic procedure to find the atoms and weights of the optimal atomic measure (see, e.g.,\, \cite{lass-book,Lau09} for details).

In addition, for the special case of the polynomial optimization problem (\ref{eq:global opt pol}), Nie \cite{Nie} shows that the flatness condition is a generic property, so that finite convergence of the lower  bounds $val^{(r)}_{outer}$ to the minimum of a polynomial over $K$ holds generically.

Note that  analogous results also hold for the stronger bounds $\overline{val}^{(r)}_{outer}$ on $val$.

\subsection{Error analysis for the case of polynomial optimization}\label{sec:convPutinar}

We now consider the special case of global polynomial optimization, i.e., problem (\ref{eq:global opt pol}), which is the case of GPM with only one affine constraint, requiring that $\mu$ is a probability measure on $K$:
$$
val=\min_{x\in K}p(x)= \min_{\mu\in\mathcal{M}(K)_+} \int_Kp(x)d\mu(x) \ \text{ s.t. } \ \int_Kd\mu(x)=1.
$$
Recall the definition of the bound $val^{(r)}_{outer}$ from (\ref{eq:upr}), which now reads
$$
val^{(r)}_{outer}= \inf_{\mu\in (\mathcal{Q}^r(g_1,\ldots,g_k))^*} \left\{\int_K p(x)d\mu(x): \int_Kd\mu(x) =1\right\}.
$$
It can be reformulated via an SDP as in (\ref{eq:upr1})-(\ref{eq:upr2}), whose dual SDP reads
\begin{equation}\label{eq:dualr}
\sup_{\lambda\in\oR} \{\lambda: p-\lambda\in \mathcal{Q}^r(g_1,\ldots,g_k)\}.
\end{equation}
By weak duality $val^{(r)}_{outer}$ is at least the optimum value of (\ref{eq:dualr}). Strong duality holds for instance if the set $K$ has a non-empty interior (since then the primal SDP is strictly feasible), or  if there is a ball constraint present in the description of the set $K$
(as shown in \cite{JH16}). Then, $val^{(r)}_{outer}$ is also given by the program (\ref{eq:dualr}), which is the case, e.g., when $K$ is a simplex, a hypercube, or a sphere.

As we saw above, the bounds $val^{(r)}_{outer}$ 
converge asymptotically to the minimum value $val$ taken by the polynomial $p$ over the set $K$ when condition (\ref{eq:ass}) holds.
We now indicate some known results on the rate of convergence of these bounds. 

For a polynomial $p=\sum_\alpha p_\alpha x^\alpha\in \oR[x]_d$, we set
$$L_p:=\max_\alpha |p_\alpha|{\alpha_1!\cdots \alpha_n!\over |\alpha|!}.$$

\begin{trailer}{Error analysis for the bounds $val^{(r)}_{outer}$} 
\begin{theorem}\cite{NS}\label{thm:NS}
Assume $K\subseteq (-1,1)^n$. There exists a constant $c>0$ (depending only on $K$) such that, for any polynomial $p$ with degree $d$, we have
$$val-val^{(r)}_{outer} \le 6d^3n^{2d} L_p
{1\over \left(\log {r\over c}\right)^{1/c}}
\quad
\text{ for all integers }\ \  r\ge c\exp\left((2d^2n^d)^c\right).
$$
\end{theorem}

\end{trailer}

Note that this result displays a very slow convergence rate, which does not reflect the good behaviour  of the bounds often observed in practice.





On the other hand, 
a sharper error analysis  holds  for the stronger bounds $\overline{val}^{(r)}_{outer}$, obtained by using the larger set $\mathcal{Q}^r(\prod_{j\in J}g_j: J\subseteq [k])$ instead of $\mathcal{Q}^r(g_1,\ldots,g_k)$.

\begin{trailer}{Error analysis for the bounds $\overline{val}^{(r)}_{outer}$}
\begin{theorem}\cite{Schw}\label{thm:Schw}
Assume $K\subseteq (-1,1)^n$. There exists a constant $c>0$ (depending only on $K$) such that, for any polynomial $p$ with degree $d$, we have
$$val-\overline{val}^{(r)}_{outer} \le  cd^4n^{2d}L_p {1\over r^{1/c}}
\quad \text{ for all integers }\ \  r\ge c d^c n^{cd}.
$$
\end{theorem}
\end{trailer}

We now recap some known sharper results for the case of polynomial optimization over special sets $K$ like the simplex, the hypercube and the sphere.
As motivation recall that this already captures well known hard combinatorial optimization problems such as the maximum independence number in a graph.

Given a graph $G=(V=[n],E)$ 
let $\alpha(G)$ denote the largest cardinality of an independent set in $G$, i.e., of a set $I\subseteq V$ that does not contain any edge of $E$.
In fact  the parameter $\alpha(G)$ can be reformulated via polynomial   optimization  over the simplex $\Delta_n$, the hypercube $[0,1]^n$, or the unit sphere $S_n$. Indeed the following results are known:
$${1\over \alpha(G)}=\min _{x\in\Delta_n} x^T(I_n+A_G)x, \quad
\alpha(G)=\max_{x\in [0,1]^n} \sum_{i\in V} x_i- \sum_{\{i,j\}\in E} x_ix_j,$$
$${2\sqrt 2\over 3\sqrt 3} \sqrt{1-{1\over \alpha (G)}} =\max_{y\in \oR^n,z\in \oR^m} \left\{2\sum_{\{i,j\} \in \overline E}z_{ij}y_iy_j: (y,z)\in S_{n+m}\right\}
$$
(see \cite{MS65,Nesterov}).
Here  $I_n$ is the identity matrix of size $n$,  $A_G$ is the adjacency matrix of $G$ (with entries $A_{ij}=A_{ji}=1$ if $\{i,j\}\in E$ and 0 otherwise),
$\overline E$ is the set of pairs of distinct elements $i,j\in V$ such that $\{i,j\}\not\in E$ and $m=|\overline E|$.


\paragraph{Error analysis for the sphere}

We first consider  the case of the sphere $K=S_n=\{x\in\oR^n: \sum_{i=1}^nx_i^2=1\}$. Then
an error analysis for the bounds $val^{(r)}_{outer}$ is known  when $p$ is a homogeneous polynomial.

First, one may reduce to the case when $p$ has even degree. Indeed, as shown in \cite{DW}, if  $p$ has odd degree $d$ then we have
$$\max\left\{p(x): \sum_{i=1}^nx_i^2=1\right\} = {d^{d/2}\over (d+1)^{(d+1)/2}} \  \max  \left\{x_{n+1}p(x): \sum_{i=1}^{n+1}x_i^2=1\right\}. $$
Another useful observation  is that,
for a homogeneous polynomial $q$ of even degree $d$, $q$ belongs to the truncated quadratic module of the sphere:
$$\mathcal{Q}^r\left(\pm\left(1-\sum_{i=1}^nx_i^2\right)\right)= \Sigma_r+\left(1-\sum_{i=1}^nx_i^2\right)\oR[x] $$ if and only if
the polynomial $q(x)\left(\sum_{i=1}^nx_i^2\right)^r$ is a sum of squares of polynomials (see \cite{KLP0}).
Therefore, when $p$ is a homogeneous polynomial of even degree $d=2a$, the parameter $val^{(r)}_{outer}$ can be reformulated as

\begin{equation}\label{psphere}
val^{(r)}_{outer}
= \min\left\{t: t\in \oR,\ t\left(\sum_{i=1}^nx_i^2\right)^r -\left(\sum_{i=1}^nx_i^2\right)^{r-a} p(x) \in \Sigma_r\right\}.
\end{equation}
Based on this,   the following error bounds for the parameters $val^{(r)}_{outer}$ are shown  in \cite{Fay,DW} (for general  polynomials) and in \cite{KLP} (for even polynomials).

\begin{theorem} 
Let   $p$ be a homogeneous polynomial of even degree $d$.
\begin{itemize}
\item[(i)] \ (\cite{Fay,DW})
There exist constants $C_{n,d}$ and $r_{n,d}$ (depending on $n$ and $d$)
such that 
$$\min_{x\in S_n}p(x)-val^{(r)}_{outer} \le {C_{n,d}\over r} \quad \text{ for all integers } r\ge r_{n,d}.
$$

\item[(ii)] \ (\cite{KLP})
If   $p$ is an even polynomial (i.e., of the form $p=\sum_{\alpha\in\oN^n_{d/2}} p_\alpha x^{2\alpha}$), then the above holds where the constant $C_{n,d}$  depends only on $d$ and  $r_{n,d}=d$.
\end{itemize}
\end{theorem}

We briefly discuss the approach in \cite{DW}, which in fact provides an error analysis for the larger range $val^{(r)}_{inner}-val^{(r)}_{outer}$.

For an integer $a$ let $\text{MSym}((\oR^n)^{\otimes a})$ denote the set of matrices acting on $(\oR^n)^{\otimes a}$ that are maximally symmetric, which means the associated $2a$-tensor is fully symmetric (i.e., invariant under the action of the symmetric group $\text{Sym}(2a)$). 
Any  homogeneous polynomial $p$ of degree $2a$ can be written as $p(x)=(x^{\otimes a})^T Z_p x^{\otimes a}$ for a (unique) $Z_p\in \text{MSym}((\oR^n)^{\otimes a})$.
Then, defining the polynomial $p_r(x)=(\sum_ix_i^2)^{r-a}p(x)$,  the program (\ref{psphere}) 
can be reformulated as 
$$val^{(r)}_{outer}= \min \left\{ \langle Z_{p_r},M\rangle: M \succeq 0, \text{Tr}(M)=1,\ M \in \text{MSym}((\oR^n)^{\otimes r})\right\}.$$
Let $M$ be an optimal solution to this program. As $M\succeq 0$ the polynomial 
$ (x^{\otimes r})^T M x^{\otimes r}$ is a sum of squares. One can scale it to obtain $h\in \Sigma_r$ which provides a probability density function on $S_n$, i.e.,
$\int_{S_n}h(x)d\mu_0(x)=1$ (with $\mu_0$ the surface measure on $S_n$), and thus $val^{(r)}_{inner} \le \int_{S_n} h(x)d\mu_0$.
Using the orthogonal polynomial basis with respect to $\mu_0$ (consisting of spherical harmonic polynomials), Doherty and Wehner \cite{DW} show a de Finetti type result, which permits to  upper bound
the range $ \int_{S_n} h(x)d\mu_0- \langle Z_{p_r},M\rangle$ and thus $ val^{(r)}_{inner}-val^{(r)}_{outer}$.

\paragraph{Error analysis for the simplex and the hypercube}

For the  simplex $K=\Delta_n=\{x\in\oR^n: x_i\ge 0\ (i\in [n]), 1-\sum_{i=1}^nx_i=0\}$ and  the hypercube $K=[0,1]^n=\{x\in \oR^n: x_i\ge 0, 1-x_i\ge 0 \ (i\in [n])\}$, a refined error analysis is known  only for the stronger bounds $\overline{val}^{(r)}_{outer}$, where we use the larger quadratic module generated by all pairwise products of the constraints defining $K$.

\begin{trailer}{Error analysis for the simplex}
\begin{theorem}\cite{KLP}
Assume $K=\Delta_n$ and $p$ is a homogeneous polynomial with degree $d$. Then we have
$$
\min_{x\in\Delta_n}p(x)-  \overline{val}^{(r)}_{outer}\le {C_d\over r} \left(\max_{x\in \Delta_n}p(x)-\min_{x\in\Delta_n}p(x)\right)\quad \text{ for all } r\ge d,
$$
where $C_d>0$ is an absolute constant depending only on $d$.
\end{theorem}
\end{trailer}

\begin{trailer}{Error analysis for the hypercube}
\begin{theorem}
\cite{KL10}
Assume $K=[0,1]^n$. 
For any polynomial $p$ with degree $d$  we have
$$
\min_{x\in [0,1]^n} p(x) - \overline{val}^{(r)}_{outer}
\le
n^d{d+1\choose 3} L_p {1\over r}
\quad \text{ for all } r\ge d.
$$
\end{theorem}
\end{trailer}

The above results show that in Theorem \ref{thm:Schw} one may choose the unknown constant to be $c=1$ (roughly) if $K$ is a hypercube or simplex.
In  both cases the proof relies on showing this error analysis for a weaker bound, which is obtained by using only nonnegative scalar multipliers (instead of sum-of-squares multipliers) in the definition of the quadratic module. See \cite{KLP,KL10} for details.

\section{Concluding remarks}
We conclude with a few remarks on  available software and future research directions.

\paragraph{Software}
The bounds based on the outer approximations (\ref{eq:upr}) described here have been implemented in the software \emph{Gloptipoly3} \cite{gloptipoly}.
The software can in fact deal with a more general version of the GPM \eqref{prob:GPM} than presented here. Namely it can deal with the problem
\[
val = \inf_{\mu_i \in \mathcal{M}(K_i)_+ \; \forall i \in \{0\}\cup [m] } \left\{ \int_{K_0} f_0(x)d\mu_0(x) \; : \; \int_{K_i} f_i(x)d\mu_i(x) = b_i \quad \forall i \in [m]\right\},
\]
where  we have a  variable measure $\mu_i\in \mathcal{M}(K_i)_+$ for each index $i\in \{0,\ldots,m\}$, with $K_i\subseteq \oR^n$ being
basic closed semi-algebraic sets defined by  (possibly different) sets of polynomial inequalities.


Due to the sizes of the resulting semidefinite programs that are solved, applicability is typically limited to $n \le 20$ variables and
low order, say $r \le 4$. This is due to the fact the matrix variables in the semidefinite programs are roughly of order ${n+r \choose r}$.
Solving larger instances requires exploiting additional structure (like sparsity) leading to more economical  semidefinite programs.
We refer, e.g.,\, to \cite{lass-book} and references therein for further details.

\paragraph{Error bounds for the inner approximation hierarchy}
The known error bounds for the inner approximation, presented earlier in Table \ref{tab:convergence rates}, are for specific choices of the set $K$
and reference measure $\mu_0 \in  \mathcal{M}(K)_+$. More work is required to understand the role of the reference measure in the convergence analysis,
and to extend the regime in $O(1/r^2)$ to more classes of sets $K$. 
In particular, an obvious choice is {whether one can sharpen the  analysis of} the the convergence rate for the Euclidean unit sphere.
As explained, such results would also have implications for grid search on cubature points on the sphere.
Cubature on the sphere is a vast research topic (see, e.g.,\, \cite[Chapter 6]{Dai-Xu book}), even in the special case of spherical $t$-designs \cite[\S6.5]{Dai-Xu book}, where all cubature weights are equal and positive.
Moreover, the complexity of polynomial optimization on spheres is not fully understood; indeed the problem is NP-hard, but allows polynomial-time approximation schemes
in special cases (see \cite{KLP,K}). 
{Sharpening the analysis of the inner approximations for polynomial optimization over spheres} 
may help to gain a more complete understanding.

\paragraph{Error bounds for the outer approximation hierarchy}
The bounds based on the outer approximation presented here are more practically suited for computation, in particular since they (sometimes) enjoy finite convergence and permit to extract the global minimizers; moreover,  as mentioned above, the dedicated software \emph{Gloptipoly3} is available
for this purpose. On the other hand, the known results on the rate of convergence are somewhat disappointing (as discussed in Section \ref{sec:convPutinar}), and in general much weaker than those known for the inner
approximation. There is certainly room for a breakthrough here; new ideas are needed to obtain convergence rates that match
the performance observed in practice.

\subsubsection*{Acknowledgement}
The authors would like to thank Fernando Mario de Oliveira Filho for insightful discussions on the duality theory of the GPM.


\begin{thebibliography}{}

\bibitem{Akh}
Akhiezer, N.I.
{\em The classical moment problem.}
Hafner, New York (1965)

\bibitem{Barvinok_course in_convexity}
Barvinok, A. {\em A course in Convexity.} Graduate Study in Mathematics, Volume 54, AMS, Providence, Rhode Island (2002)

\bibitem{BT}
Bayer, C., and Teichmann, J.
The proof of Tchakaloff's theorem.
{\em Proceedings of
the American Mathematical Society},
134:3035--3040 (2006)

\bibitem{BenTal}
Ben Tal, A., and Nemirovski, A.
{\em Lectures on Modern Convex Optimization: Analysis, Algorithms, and Engineering Applications.}
{ MPS/SIAM Series on Optimization}, 2, SIAM (2001)

\bibitem{BGGLT}
Bhattiprolu, V., Ghosh, M., Guruswami, V., Lee, E., Tulsiani, M.
Weak decoupling, polynomial folds and approximate optimization over the sphere.
arXiv:1611.05998v2 (2017)


\bibitem{Castella}
Castella, M. Rational optimization for nonlinear reconstruction
with approximate
$\ell_0$ penalization. Preprint version available at	arXiv:1808.00724 (2018)



\bibitem{Cools Encyclopedia cubature}
 Cools, R. An encyclopaedia of cubature formulas, \emph{J. Complexity}, 19, 445--453 (2003).

\bibitem{CF}
Curto, R.E., and  Fialkow, L.A.
Solution of the truncated complex moment problem for flat data,
{\em Memoirs
of the American Mathematical Society}
119
(568) (1996)


\bibitem{Dai-Xu book}
Dai, F., and Xu, Y.
\emph{Approximation Theory
and Harmonic Analysis
on Spheres and Balls}, Springer, New York (2013)

\bibitem{DHL SIOPT}
{De Klerk}, E., Hess, R., and Laurent, M.
Improved convergence rates for Lasserre-type hierarchies of upper bounds for box-constrained polynomial optimization.
{\em SIAM Journal on Optimization}
 27(1), 347--367 (2017)

\bibitem{K}
De Klerk, E.
The complexity of optimizing over a simplex, hypercube or sphere: a short survey.
{\em Central European Journal of Operations Research}, 16(2), 111--125 (2008)


\bibitem{KLLS17}
De Klerk, E.,  Lasserre, J.B.,  Laurent, M., and  Sun, S.
Bound-constrained polynomial optimization using only elementary calculations.
{\em Mathematics of Operations Research,} 42(3), 834--853 (2017)


\bibitem{KL10}
{De Klerk}, E., Laurent, M.
Error bounds for some semidefinite programming approaches to polynomial minimization on the hypercube.
{\em  SIAM Journal on Optimization} 20(6), 3104--3120 (2010)

\bibitem{KL_MOR_2017}
{De Klerk}, E., Laurent, M. Comparison of Lasserre's measure-based  bounds for polynomial optimization to bounds obtained by simulated annealing. \emph{Mathematics of Operations Research}, to appear. Preprint version available at \url{http://arxiv.org/abs/1703.00744} (2017)

\bibitem{KL_2018}
{De Klerk}, E., Laurent, M. Worst-case examples for Lasserre's measure--based hierarchy for polynomial optimization on the hypercube.  	
 \emph{Mathematics of Operations Research}, to appear. Preprint version available at \url{http://arxiv.org/abs/1804.05524} (2018)

\bibitem{KLP0}
{De Klerk}, E., Laurent, M., and Parrilo, P.
On the equivalence of algebraic approaches to the miniization of forms on the simplex.
In D. Henrion and A. Garulli (eds), {\em Positive Polynomials in Control}, 121--133, Springer (2005)


\bibitem{KLP}
{De Klerk}, E., Laurent, M., and Parrilo, P.
A PTAS for the minimization of polynomials of fixed degree over the simplex.
{\em Theoretical Computer Science,} 361(2-3), 210--225 (2006)

\bibitem{KLS_MPA}
{De Klerk}, E., Laurent, M., Sun, Z. Convergence analysis for Lasserre's measure-based hierarchy of upper bounds for polynomial optimization. \textit{Mathematical Programming Series A} 162(1), 363--392 (2017)

\bibitem{KLS15}
{De Klerk}, E., Laurent, M., Sun, Z.
An error analysis for polynomial optimization over the simplex based on the multivariate hypergeometric distribution. {\em SIAM Journal on Optimization,} 25(3), 1498--1514 (2015)


\bibitem{KLSV}
De Klerk, E., Laurent, M., Sun, Z., and Vera, J.
On the convergence rate of grid search for polynomial optimization over the simplex.
{\em Optimization Letters,} 11(3), 597--608 (2017)


\bibitem{DK-P-K dist robust opt}
De Klerk, E., Postek, K., and Kuhn, D.
Distributionally robust optimization with polynomial densities: theory, models and algorithms.
Preprint version available at	arXiv:1805.03588 (2018)

\bibitem{DW}
Doherty, A.C., Wehner, S.
Convergence of SDP hierarchies for polynomial optimization on the hypersphere.
arXiv:1210.5048v2 (2013)

\bibitem{Dunkl}
 Dunkl, C.F.,  and  Xu., Y.
{\em
Orthogonal Polynomials of Several Variables,}
Cambridge University Press (2001)

\bibitem{Fay}
Faybusovich, L.
Global optimization of homogeneous polynomials on the simplex and on the sphere.
In C. Floudas and P. Pardalos (eds), {\em Frontiers in Global Optimization}, Kluwer (2003)

\bibitem{Folland}
Folland, G.B.
How to integrate a polynomial over a sphere?
{\em The American Mathematical Monthly}, 108(5), 446--448 (2001)

\bibitem{Han}
Handelman, D.
Representing polynomials by positive linear functions on compact convex polyhedra.
{\em Pacific J. Math.}, 132(1), 35--62 (1988)

\bibitem{gloptipoly}
 Henrion, D.,  Lasserre, J.B., and Loefberg, J. GloptiPoly 3: moments, optimization and semidefinite programming.
  \emph{Optimization Methods and Software}, 24(4-5), 761--779 (2009). Software download: \url{www.laas.fr/$\sim$henrion/software/gloptipoly3}

\bibitem{JH16}
Josz, C.,   Henrion, D.
Strong duality in Lasserre's hierarchy for polynomial optimization.
{\em Optim. Letters}, 10, 3--10 (2016)

\bibitem{Kalai-Vempala 2006}
 Kalai. A.~T., and Vempala, S.
 Simulated annealing for convex optimization.
 \emph{Mathematics of Operations Research}, 31(2),
  253--266 (2006)

\bibitem{Kemperman}
Kemperman, J.H.B.
The general moment problem, a geometric approach.
{\em The Annals of Mathematics Statistics}, 39, 93--122 (1968)

\bibitem{Landau}
Landau, H.
{\em Moments in Mathematics}, {\em Proc. Sympos. Appl. Math.}, 37 (1987)

\bibitem{lass-siopt-01}
Lasserre, J.B.
Global optimization with polynomials and the problem of moments. \emph{SIAM J. Optim.} 11, 796--817 (2001)

\bibitem{Las2008}
Lasserre, J.B. A semidefinite programming approach to the generalized problem of moments.
\emph{Mathematical Programming Series B} 112, 65--92 (2008)

\bibitem{lass-book}
Lasserre, J.B.
{\em Moments, Positive Polynomials and Their Applications.} Imperial College Press (2009)

\bibitem{lass-book2}
Lasserre, J.B.
{\em Introduction to Polynomial and Semi-Algebraic Optimization}. Cambridge University Press  (2015)

\bibitem{Las11}
Lasserre, J.B. A new look at nonnegativity on closed sets and polynomial optimization. {\em SIAM Journal on Optimization} 21(3), 864--885 (2011)

\bibitem{Las18}
Lasserre, J.B.
The moment-SOS hierarchy.
{\em Proc. Int. Cong. of Math. ? 2018},
Rio de Janeiro, 3, 3761--3784 (2018)

\bibitem{Lau09}
Laurent, M.
Sums of squares, moment matrices and optimization over polynomials. In {\em Emerging Applications of Algebraic Geometry,} Vol. 149 of IMA Volumes in Mathematics and its Applications, M. Putinar and S. Sullivant (eds.), Springer, 157-270 (2009)

\bibitem{Lau14}
Laurent, M.
Optimization over polynomials: Selected topics. In Chapter 16 (Control Theory and Optimization) of {\em Proc. Int. Cong. of Math. 2014}.  Jang, S. Y., Kim, Y. R., Lee, D-W. \& Yie, I. (eds.). Seoul: Kyung Moon SA Co. Ltd., p. 843--869 (2014)

\bibitem{Martinez et al quadrature}
Martinez, A.,  Piazzon, F.,  Sommariva, A.,  and Vianello, M.
Quadrature-based polynomial
optimization. Manuscript (2018). \url{http://www.math.unipd.it/~marcov/pdf/quadropt.pdf}

\bibitem{MS65}
Motzkin, T.S., Sraus, E.G.
Maxima for graphs and a new proof of a theorem of T\'uran.
{\em Canadian J. Math.}, 17, 533--540 (1965)

\bibitem{Nesterov}
Nesterov, Yu.
Random walk in a simplex and quadratic optimization over convex polytopes.
CORE Discussion Paper 2003/71, CORE-UCL, Louvain-La-Neuve (2003)

\bibitem{Nie}
Nie, J.
 Optimality conditions and finite convergence of Lasserre's hierarchy.
{\em Mathematical Programming, Ser. A,}146(1-2), 97--121 (2014)

\bibitem{NS}
Nie, J.,  and Schweighofer, M.
On the complexity of Putinar's positivstellensatz
{\em Journal of Complexity} 23, 135--150 (2007)

\bibitem{PV}
Piazzon, F., Vianello, M.
Markov inequalities, Dubiner distance,
norming meshes and polynomial
optimization on convex bodies.
Preprint at Optimization Online (2018)

\bibitem{PR}
Powers, V., Reznick, B.
A new bound for P\'olya's theorem with applications to polynomials positive on polyhedra.
{\em J. Pure Appl. Algebra}, 164, 221--229 (2001)

\bibitem{Putinar}
Putinar, M.
Positive polynomials on compact semi-algebraic sets.
{\em Ind. Univ. Math. J.} 42, 969--984 (1993)

\bibitem{Putinar cubature}
Putinar, M.
A note on Tchakaloff's theorem.
Proceedings of the
American Mathematical Society,
 125(8), 2409--2414 (1997)

\bibitem{rogosinski}
Rogosinski, W.W. Moments of non-negative mass, \textit{Proceedings of the Royal Society A} 245, 1--27 (1958)

\bibitem{Ryu_Boyd_quadrature}
 Ryu, E.K. and Boyd, S.P.
Extensions of Gauss Quadrature Via Linear
Programming. \emph{Foundations of Computational Mathematics}
15(4),  953--971 (2015)

\bibitem{Schm}
Schm\"udgen, K.
The $K$-moment problem for compact semi-algebraic sets.
{\em Math. Ann.}, 289, 203--206 (1991)


\bibitem{Sch}
Schm\"udgen, K.
{\em The moment problem}.
Springer (2017)

\bibitem{Schw}
Schweighofer, M.
On the complexity of Schmüdgen's Positivstellensatz
{\em Journal of Complexity} 20(4), 529-543 (2004)



\bibitem{Thompson's problem}
Schwartz, R.E.
The 5 electron case of Thomson's problem. \emph{Exp. Math} 22(2), 157--186 (2013)


\bibitem{Shapiro2001}
Shapiro, A. \emph{On duality theory of conic linear problems}, Semi-Infinite Programming: Recent Advances (M,{\'A}. Goberna and M.A. L{\'o}pez,
  eds.), Springer, 135--165 (2001)

\bibitem{Tao_epsilon}
Tao, T. An Epsilon of Room, I: Real Analysis: pages from year three of a mathematical blog. AMS, Graduate Studies in Mathematics
Volume: 117 (2010)

\bibitem{Trefethen cubature review}
Trefethen, L.N.
Cubature, Approximation, and
Isotropy in the Hypercube. \emph{SIAM Review}, 59(3), 469--491 (2017)

\bibitem{Tchakaloff}
Tchakaloff, V. Formules de cubature \'{m}ecanique \`{c}oefficients non n\'{e}gatifs, \emph{Bull. Sci. Math.},
81, 123--134 (1957)

\end{thebibliography}
\end{document}